\documentclass{article}
\usepackage{amsthm}
\usepackage{amsmath}
\usepackage{amssymb}
\usepackage{tikz}
\usepackage{booktabs}
\usepackage{setspace} 
\usepackage{enumerate}
\usepackage{geometry}
\usepackage{graphicx}
\usepackage{xcolor}
\usepackage{ifthen}
\usepackage[hidelinks]{hyperref}

\usepackage{float}
\floatstyle{ruled}
\newfloat{algorithm}{h}{loa}
\floatname{algorithm}{Algorithm}
\usepackage{algpseudocode}
\usepackage{tabularx}
\usepackage{booktabs}
\usepackage{array}

\newtheorem{prop}{Proposition}
\newtheorem{rem}[prop]{Remark}

\newtheorem{lem}[prop]{Lemma}

\DeclareMathOperator{\TV}{TV}

\DeclareMathOperator{\TGV}{TGV}

\DeclareMathOperator*{\argmin}{argmin}

\newcommand{\TGVat}{\TGV_\alpha ^2}

\newcommand{\ctc}{\emph{CT-cvx}}
\newcommand{\ctsc}{\emph{CT-scvx}}
\newcommand{\txt}{\emph{TXT}}
\newcommand{\tgv}{\emph{TGV}}
\newcommand{\bmtd}{\emph{BM3D}}
\newcommand{\cl}{\emph{CL}}

\DeclareMathOperator{\range}{Rg}

\DeclareMathOperator{\proj}{proj}

\DeclareMathOperator{\prox}{prox}
\DeclareMathOperator{\diag}{diag}

\newcommand{\reshape}[2][NM,nn]{[#2_{(#1)}]}
\newcommand{\B}{ \mathcal{B} }

\newcommand{\st}{ \, \left|\right.\, }

\newcommand{\rank}{\text{rank}}

\newcommand{\tp}{\otimes_\pi}
\newcommand{\itp}{\otimes_\mathfrak{i}}
\newcommand{\itn}{\mathfrak{i}}

\newcommand{\ctp}{\otimes_\pi}

\makeatletter
\newcommand{\vc}[2][r]{%
  \gdef\@VORNE{1}
  \left(\hskip-\arraycolsep%
    \begin{array}{#1}\vekSp@lten{#2}\end{array}%
  \hskip-\arraycolsep\right)}
\def\vekSp@lten#1{\xvekSp@lten#1;vekL@stLine;}
\def\vekL@stLine{vekL@stLine}
\def\xvekSp@lten#1;{\def\temp{#1}%
  \ifx\temp\vekL@stLine
  \else
    \ifnum\@VORNE=1\gdef\@VORNE{0}
    \else\@arraycr\fi%
    #1%
    \expandafter\xvekSp@lten
  \fi}
\makeatother

\newcommand{\Nuc}{\mathcal{N}}
\newcommand{\Lin}{\mathcal{L}}

\newcommand{\nuc}{\text{nuc}}

\newcommand{\ncl}{N_\text{cl,s}}

\newcommand{\R}{\mathbb{R}}
\newcommand{\N}{\mathbb{N}}

\newcommand{\symgrad}{\mathcal{E}}
\newcommand{\wrt}{\:\mathrm{d}}

\newcommand{\beq}{\begin{equation}}
\newcommand{\eeq}{\end{equation}}

\newcommand{\I}{\mathcal{I}}
\newcommand{\M}{\mathcal{M}}

\title{A convex variational model for learning convolutional image atoms from incomplete data}

\author{A. Chambolle, M. Holler and T. Pock}
\date{}

\begin{document}
\maketitle
\begin{abstract}
A variational model for learning convolutional image atoms from corrupted and/or incomplete data is introduced and analyzed both in function space and numerically.
Building on lifting and relaxation strategies, the proposed approach is convex and allows for simultaneous image reconstruction and atom-learning in a general, inverse problems context. 
Further, motivated by an improved numerical performance, also a semi-convex variant is included in the analysis and the experiments of the paper. 
For both settings, fundamental analytical properties allowing in particular to ensure well-posedness and stability results for inverse problems are proven in a continuous setting. 
Exploiting convexity, globally optimal solutions are further computed numerically for applications with incomplete, noisy and blurry data and numerical results are shown.
\end{abstract}


\section{Introduction}

An important task in image processing is to achieve an appropriate regularization or smoothing of images or image related data. In particular, this is indispensable for most application-driven problems in the field, such as denoising, inpainting, reconstruction, segmentation, registration or classification. Also beyond imaging, for general problem settings in the field of inverse problems, an appropriate regularization of unknowns plays a central role as it allows for a stable inversion procedure.

Variational methods and partial-differential-equation-based methods can now be regarded as classical regularization approaches of mathematical image processing (see for instance  \cite{Weickert98anisotropic,aubert2006mathematical,scherzer2009variationalmethods,mallat2009wavelettour}). An advantage of such methods is the existence of a well-established mathematical theory and, in particular for variational methods, a direct applicability to general inverse problems with provable stability and recovery guarantees \cite{Hofmann07_nonlinear_tikhonov_banach,holler17coupled_reg}.
While in particular piecewise smooth images are typically well-described by such methods, their performance for oscillatory or texture-like structures, however, is often limited to pre-described patterns (see for instance \cite{Bredies2018oscillation}).

Data-adaptive methods such as patch-based methods (see for instance \cite{Lebrun12_image_denoising_review,dabov2007bm3d,Buades05,Delon17_mixing_tv_non_local}) on the other hand are able to exploit redundant structures in images independent of an a-priory description and are, at least for some specific tasks, often superior to variational- and PDE-based methods. In particular machine-learning-based methods have advanced the state-of-the art significantly in many typical imaging applications in the past years. A disadvantage of such methods, however, is a current lack of mathematical understanding in particular compared to existing results for variational methods in the context of inverse problems. Data-adaptive methods typically build on learning some image atoms in some way or the other and neither the mapping of training data to image atoms nor the application of learned atoms to reconstruction is known to enjoy the same stability and regularization properties as with classical methods. In particular, for typical deep neuronal networks, the training step corresponds to the minimization of a non-convex energy, but neither is this minimization carried out until convergence nor are the properties of local minimizers well understood mathematically. 

To goal of this work is to provide a first step towards bridging the gap between data-driven methods and variational methods. Building on a tensorial-lifting approach, we introduce a convex variational method for learning image atoms from noisy and/or incomplete data in an inverse problems context. We further extend this model by a semi-convex variant that improves the performance in some applications. For both settings, we are able to prove well-posedness results in function space and, for the convex version, to compute globally optimal solutions numerically. In particular, classical stability and convergence results for inverse problems such as the ones of \cite{Hofmann07_nonlinear_tikhonov_banach,holler17coupled_reg} are applicable to our model, providing a stable recovery of both learned atoms and images from given, incomplete data.

Our approach is motivated by a sparse, convolutional representation of images via a few image atoms and conceptually allows for a joint learning of image atoms and image reconstruction in a single step. Nevertheless, it can also be regarded purely as an a approach for learning image atoms from potentially incomplete data in a training step, after which the learned atoms can be further incorporated in a second step, e.g., for reconstruction or classification.
It should also be noted that, while we show some examples where our approach achieves a good practical performance for image reconstruction compared to existing methods, the main purpose of this paper is to provide a mathematical understanding rather than an algorithm that achieves the best performance in practice. 

Regarding existing literature in the context of data-adaptive variational learning approaches in imaging, we note that there are many recent approaches that aim to infer either parameter or filters for variational methods from given training data, see e.g. \cite{Kunisch13bilevel,Calatroni16bilevel,Hintermueller17_bileveltv}. 
A continuation of such techniques more towards the architecture of neuronal networks are so called variational networks \cite{kobler2017variational,adler2018learned} where not only model parameters but also components of the solution algorithm such as stepsizes or proximal mappings are learned. We also refer to \cite{lunz2018adversarial} for a recent work on combining variational methods and neuronal networks.
While for some of those methods also a function space theory is available, the learning step is still non-convex and one can in general only expect to obtain locally optimal solutions.

\subsection{Outline of the paper}

In Section \ref{sec:main_approach} we present the main ideas for our approach in a formal setting. This is done from two perspectives, once from the perspective of a convolutional lasso approach and once from the perspective of patch-based methods. In Section \ref{sec:continuous_setting} we then carry out an analysis of the proposed model in function space, were we motivate our approach via convex relaxation and derive well-posedness results. Section \ref{sec:discrete_setting} then presents the model in a discrete setting and the numerical solution strategy and Section \ref{sec:numerical_results} provides numerical results and a comparison to existing methods. At last, an appendix provides a brief overview on some results for tensor spaces that are used in Section \ref{sec:continuous_setting}. We note that, while the analysis of Section 3 is an important part of our work, the paper is structured in a way such that readers only interested in the conceptual idea and the realization of our approach can skip Section \ref{sec:continuous_setting} and focus on Sections \ref{sec:main_approach} and \ref{sec:discrete_setting}.

\begin{figure}
\centering
\scalebox{0.95}{
\begin{tikzpicture}
%Base image
\draw[blue,thick] (0,0) --node[below,black]{Image} (5.5,0);
\node[below] at (2.75,3.5)    {$u=KC$};

%Selected patches
\foreach \x in {-0.5,1.0,2.25,4.5}
{
	\draw[dotted] (\x,2.5-0.3*\x) -- (\x,0.2);
	\draw[dotted] (\x+1,2.5-0.3*\x) -- (\x+1,0.2);
	\draw[green,thick] (\x,0.2) -- (\x+1,0.2);
}
\node[left] at (-1.1, 0.2)    {Active atoms};
\draw (-1,0.05) -- (-1.1,0.05) -- (-1.1,0.35) -- (-1,0.35);
%Lifted patches
\foreach \x in {-0.75,-0.5,...,5.25}
{
	\draw[red,line width=0.4pt] (\x,2.5-0.3*\x) -- (\x+1,2.5-0.3*\x);
	}
%Re-drawing of selected patches
\foreach \x in {-0.5,1.0,2.25,4.5}
{
	\draw[green,thick] (\x,2.5-0.3*\x) -- (\x+1,2.5-0.3*\x);
	}
\node[left] at (-1.1, 2.5-0.3*2.25)    {Lifted atoms};
\draw (-1,2.5+0.3*0.75+0.15) -- (-1.1,2.5+0.3*0.75+0.15) -- (-1.1,2.5-0.3*5.25-0.15) -- (-1,2.5-0.3*5.25-0.15);
\draw (6.5,2.5+0.3*0.75+0.15) -- (6.6,2.5+0.3*0.75+0.15) -- (6.6,2.5-0.3*5.25-0.15) -- (6.5,2.5-0.3*5.25-0.15);
%Patch-matrix
\foreach \x in {-0.75,-0.5,...,5.25}
{
	\draw[red,line width=0.4pt] (8,2.5-0.3*\x) -- (9,2.5-0.3*\x);
	}
%Re-drawing of selected patches
\foreach \x in {-0.5,1.0,2.25,4.5}
{
	\draw[green] (8,2.5-0.3*\x) -- (9,2.5-0.3*\x);
	}
\node[below] at (8.5,3.5)    {$C$};
\draw[dotted] (7.85,0.8) -- (9.15,0.8) -- (9.15,2.85) -- (7.85,2.85) -- (7.85,0.8);
\node[below] at (9.8,0)    {Atom-matrix and decomposition};
%Arrow
\draw[<-] (7,2.5-0.3*2.25) -- (7.5,2.5-0.3*2.25);	
%Equality and product
\node at (9.5,2.5-0.3*2.25) {$=$};	
\node at (10.5,2.5-0.3*2.25) {$\otimes$};	
\node[below] at (10.5,3.5) {$c \otimes p$};	
%Coefficient image and re-drawing
\foreach \x in {-0.75,-0.5,...,5.25}
{
	\draw[red,line width=0.4pt] (10,2.5-0.3*\x-0.0) rectangle (10.1,2.5-0.3*\x + 0.0);
}
\foreach \x in {-0.5,1.0,2.25,4.5}
{
	\draw[green,fill] (10,2.5-0.3*\x-0.01) rectangle (10.1,2.5-0.3*\x + 0.01);
}
\draw[dotted] (9.85,0.8) -- (10.25,0.8) -- (10.25,2.85) -- (9.85,2.85) -- (9.85,0.8);
\draw[green,thick] (10.9,2.5-0.3*2.25) -- (11.9,2.5-0.3*2.25);
\draw[dotted] (10.8,2.5-0.3*2.25+0.1) rectangle (12,2.5-0.3*2.25-0.1);
\end{tikzpicture}
}
\caption{\label{fig:lifting_visualization}Visualization of the atom-lifting approach for 1D images. The green  (thick) lines in the atom-matrix correspond to non-zero (active) atoms and are placed in the image at the corresponding positions.}
\end{figure}
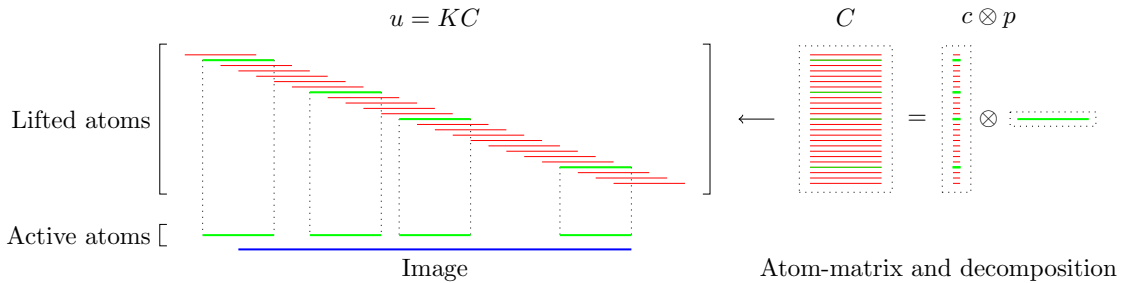

\section{A convex approach to image atoms} \label{sec:main_approach}
In this section we present the proposed approach to image-atom-learning and texture reconstruction, where we focus on explaining the main ideas rather than precise definitions of the involved terms. For the latter, we refer to Section \ref{sec:continuous_setting} for the continuous model and Section \ref{sec:discrete_setting} for the discrete setting.

Our starting point is the convolutional lasso problem \cite{Zeiler10convlasso,Chambolle16acta}, which aims to decompose a given image $u$ as a sparse linear combination of basic atoms $(p_i)_{i=1}^k$ with coefficient images $(c_i)_{i=1}^k$ by inverting a sum of convolutions as follows
\[ \min_{ (c_i)_i, (p_i)_i }  \sum_{i=1}^k \|c_i\|_1   \qquad \text{s.t. }
\begin{cases}
u = \sum_{i=1}^k c_i * p_i, \\
\|p_i\| _2 \leq 1 \text{ for } i=1,\ldots,k.
\end{cases}
\]
It is important to note that, by choosing the $(c_i)_i$ to be composed of delta peaks, this allows to place the atoms $(p_i)_i$ and any position in the image.
In \cite{Zeiler10convlasso}, this model has been used in the context of convolutional neural networks for generating image atoms and other image related tasks. Subsequently, many works have dealt with the algorithmic solution of the resulting optimization problem, where the main difficulty lies in the non-convexity of the atom learning step, and we refer to \cite{Wohlberg18conv_dict_learning_review} for a recent review.

Our goal is to obtain a convex relaxation of this model that can be used for both, learning image atoms from potentially noisy data as well as image reconstruction tasks such as inpainting, deblurring or denoising. 
To this aim, we lift the model to the tensor product space of coefficient images and image atoms, i.e., the space of all tensors $C = \sum_{i} c_i \otimes p_i$ with $c_i \otimes p_i$ being a rank-1 tensor such that $(c_i \otimes p_i)(x,y) = c_i(x)p_i(y)$. We refer to Figure \ref{fig:lifting_visualization} for a visualization of this lifting in a one-dimensional setting, where both coefficients and image atoms are vectors and $c_i \otimes p_i$ corresponds to a rank-one matrix. Notably, in this tensor product space, the convolution $c_i *p_i$ can be written as linear operator $K$ such that $KC(x) = \sum_{i} K(c_i \otimes p_i)(x) =  \sum_i \int p_i(x-y)c_i(y)$. Exploiting redundancies in the penalization of $(\|c_i\|_1)_i $ and the constraint $\|p_i\|_2 \leq 1$, $i=1\ldots,k$ and re-writing the above setting in the lifted tensor space, as discussed in Section \ref{sec:continuous_setting}, we obtain the following minimization problem as convex relaxation of the convolutional lasso approach
\[ \min_{ C }   \|C\|_{1,2} \qquad \text{s.t. }u =  KC,\]
where $\|\cdot \|_{1,2}$ takes the $1$-norm and $2$-norm of $C$ in coefficient and atom direction, respectively. Now while a main feature of the original model was that the number of image atoms was fixed, this is no longer the case in the convex relaxation and would correspond to constraining the rank of the lifted variable $C$ (defined as the minimal number of simple tensors needed to decompose $C$) to be below a fixed number. As convex surrogate, we add an additional penalization of the nuclear norm of $C$ in the above objective functional (here we refer to the nuclear norm of $C$ in the tensor product space which, in the discretization of our setting, coincides with the classical nuclear norm of a matrix-reshaping of $C$). Allowing also for additional linear constraints on $C$ via a linear operator $M$, we arrive at the following convex norm that measures the decomposability of a given image $u$ into a sparse combination of atoms as
\[ N_\nu(u)  = \min _{C}\, \nu \|C\|_{1,2 } + (1-\nu)\|C\|_* \qquad \text{ such that } u= KC, \, MC = 0. \]
Interestingly, this provides a convex model for learning image atoms, which for simple images admitting a sparse representation seems quite effective. In addition, this can in principle also be used as a prior for image reconstruction tasks in the context of inverse problems via solving for example
\[ \min _u \frac{\lambda}{2}\|Au-u_0\|_2^2 + N_\nu(u) ,\]
with $u_0$ given some corrupted data, $A$ a forward operator and $\lambda>0$ a parameter.

Both the original motivation for our model as well as its convex variant have many similarities with dictionary learning and patch-based methods. The next section strives to clarify similarities and difference and provides a rather interesting, different perspective on our model.

\subsection{A dictionary-learning/patch-based methods perspective}

In classical dictionary-learning-based approaches, the aim is to represent a resorted matrix of image patches as a sparse combination of dictionary atoms. That is, with $u \in \R^{NM}$ a vectorized version of an image and $D = (D_1,\ldots,D_l)^T \in \R^{l \times nm}$ a patch-matrix containing $l$ vectorized (typically overlapping) images patches of size $nm$, the goal is to obtain a decomposition $D = cd$, where $c \in \R^{l \times k}$ is a coefficient matrix and $p \in \R^{k \times nm}$ is a matrix of $k$ dictionary atoms such that $c_{i,j}$ is the coefficient for the atom $p_{j,\cdot}$ in the representation of the patch $D_{i}$. In order to achieve a decomposition in this form, using only a sparse representation of dictionary atoms, a classical approach is to solve
\[ \min _{c,p} \frac{\lambda}{2}\| cp - D \|_2^2 + \|c\|_1 + R(p), \]
where $R$ potentially puts additional constraints or cost on the dictionary atoms, e.g. $R(p) = 0$ if $\|p_{j,\cdot} \|_2 \leq 1$ for all $j $ and $R(p) = \infty $ else.

A difficulty with such an approach is again the bilinear and hence non-convex nature of the optimization problem, leading to potentially many non-optimal stationary points and making the approach sensitive to initialization.

As a remedy, one strategy is to consider a convex variant (see for instance \cite{bach2008convex}). That is, re-writing the above minimization problem (and again using the ambiguity in the product $cd$ to eliminate the  $L^2 $ constraint) we arrive at the problem
\[ \min _{C: \rank (C) \leq k} \frac{\lambda}{2}\| C - D \|_2^2 + \|C\|_{1,2}, \]
where $\|C \|_{1,2} = \sum_{i} \|C_{i,\cdot} \|_2 $. A possible convexification is then given as 
\begin{equation} \label{eq:dict_denoise_convex}
 \min _{C} \frac{\lambda}{2}\| C - D \|_2^2 + \nu\|C\|_{1,2} + (1-\nu)\| C\|_{*}, 
 \end{equation}
where $\|\cdot \|_{*}$ is the nuclear norm of the matrix $C$.

A disadvantage of such an approach is that the selection of patches is a-priory fixed and that the lifted matrix $C$ has to approximate each patch. In the situation of overlapping patches, this means that individual rows of $C$ have to represent different shifted version of the same patch several times, which inherently contradicts the low-rank assumption. 

\begin{figure}[t]
\center 

\includegraphics[width=0.25\linewidth]{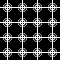}
\includegraphics[width=0.25\linewidth]{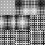}
\includegraphics[width=0.25\linewidth,trim=0 150 0 200,clip]{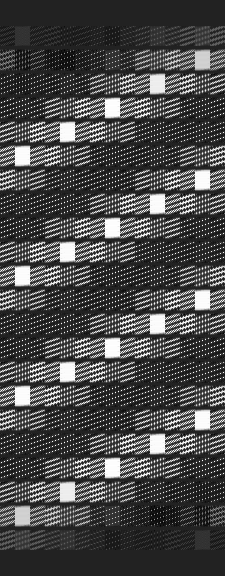}

\hspace*{0.258\linewidth}
\includegraphics[width=0.25\linewidth]{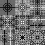}
\includegraphics[width=0.25\linewidth,trim=0 150 0 200,clip]{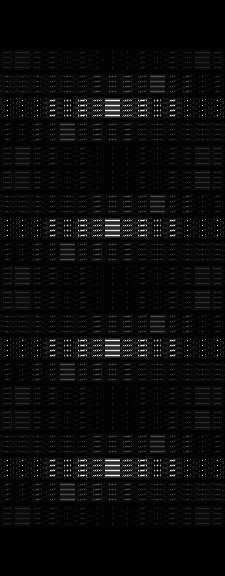}

\caption{\label{fig:ksvd_vs_patch} Patch-based representation of test images. Left: Original image, Middle: Nine most important patches for each method (top: patch denoising, bottom: patch reconstruction), Right: Section of the corresponding patch matrices.}

\end{figure}

It is now interesting to see our proposed approach in relation to this dictionary learning methods and the above-described disadvantage: Denote again by $K$ the lifted version of the convolution operator, which in the discrete setting takes a lifted patch-matrix as input and provides an image composted of overlapping patches as output. It is than easy to see that $K^*$, the adjoint of $K$, is in fact a patch selection operator and it holds that $K K^* = I$. Now using $K^*$, the approach in \eqref{eq:dict_denoise_convex} can be re-written as
\begin{equation} \label{eq:patch_denoising}
 \min _{C} \frac{\lambda}{2}\| C - K^*u \|_2^2 + \nu \|C\|_{1,2} + (1-\nu)\| C\|_{*}, 
 \end{equation}
where we remember that $u$ is the original image. Considering the problem of finding an optimal patch-based representation of an image as the problem if inverting $K$, we can see that the previous approach in fact first applies a right-inverse of $K$ and then decomposes the data. Taking this perspective, however, it seems much more natural to consider instead an adjoint formulation as
\begin{equation} \label{eq:patch_reconstruction}
 \min _{C} \frac{\lambda}{2}\| KC - u \|_2^2 + \nu \|C\|_{1,2} + (1-\nu)\| C\|_{*}. 
 \end{equation}
Indeed, this means that we do not fix the patch-decomposition of the image a-priory but rather allow the method itself to optimally select the position and size of patches. In particular, given a particular patch at an arbitrary location, this patch can be represented by using only one line of $C$ and the other lines (corresponding to shifted versions) can be left empty. Figure \ref{fig:ksvd_vs_patch} shows the resulting improvement by solving both of the above optimization problems for a particular test image, where the parameters are set such that the data error of both methods, i.e., $\|KC -U \|_2 ^2 $ is approximately the same. As can be seen, solving \eqref{eq:patch_denoising}, which we call patch denoising, does not yield meaningful dictionary atoms as the dictionary elements need to represent different, shifted version of the single patch that makes up the image. In contrast to that, solving \eqref{eq:patch_reconstruction}, which we call patch reconstruction, allows to identify the underlying patch of the image and the corresponding patch-matrix is indeed row-sparse.

\subsection{The variational model}
Now while the proposed model can, in principle, describe any kind of images, in particular its convex relaxation seems best suited for situations where the image can be described by only a few, repeating atoms, as would be for instance the case with texture images. In particular, since we do not include rotations in the model, there are many simple situations, such as $u$ being the characteristic function of a disk, which would in fact require an infinite number of atoms. To overcome this, it seems beneficial to include an additional term which is rotationally invariant and takes care of piecewise smooth parts of the image. Denoting $R$ to be any regularization functional for piecewise smooth data and taking the infimal convolution of this functional with our atom-based norm, we then arrive at the convex model
\[ \min _{u,v}  \frac{\lambda}{2}\|Au-u_0\|_2^2 + \mu_1R(u-v) + \mu_2N_\nu(v),\]
for learning image data and image atoms kernels from potentially noisy or incomplete measurements.

A particular example of this model can be given when choosing $R=\TV$, the total variation function \cite{rudin1992tv}. In this setting, a natural choice for the operator $M$ in the definition of $N_\nu$ is to take the point-wise mean of the lifted variable in atom direction, which corresponds to constraining the learned atoms to have zero mean and enforces some kind of orthogonality between the cartoon and the texture part in the spirit of \cite{meyer2001oscillating}. In our numerical experiments, in order to obtain an even richer model for the cartoon part, we use the second order Total Generalized Variation function ($\TGVat)$ \cite{bredies2010tgv,holler14inversetgv} as cartoon prior and, in the spirit of a dual $\TGVat$ norm, use $M$ to constrain the $0$th and $1$st moments of the atoms to be zero.

We also remark that, as shown in the analysis part of Section \ref{sec:continuous_setting}, while an $\ell^1/\ell^2$-type norm on the lifted variables indeed arises as convex relaxation of the convolutional lasso approach, the addition of the nuclear norm is to some extend arbitrary and in fact in the context of compressed sensing it is known that a summation of two norms is is suboptimal for a joint penalization of sparsity and rank \cite{Oymak15low_rank_sparse} (we refer to Remark \ref{rem:nuclear_norm_hilbert} below for an extended discussion). Indeed, our numerical experiments also indicate that the performance of our method is to some extend limited by a sub-optimal relaxation of a joint sparsity and rank penalization. To account for that, we also tested with semi-convex potential functions (instead of the identity) for a penalization of the singular values in the nuclear norm. Since this provided a significant improvement in some situations, we also include this more general setting in our analysis and the numerical results.

\section{The model in a continuous setting} \label{sec:continuous_setting}

The goal of this section is to define and analyze the model introduced in Section \ref{sec:main_approach} in a continuous setting.
To this aim, we regard images as functions in the Lebesgue spaces $L^q(\Omega)$ with a bounded Lipschitz domain $\Omega \subset \R^d$, $d \in \N$ and $1 < q \leq 2$. Image atoms are regarded as functions in $L^s(\Sigma)$, with $\Sigma \subset \R^d$ a second (smaller) bounded Lipschitz domain (either a circle or a square around the origin) and $s \in [q,\infty]$ an exponent that is a-priori allowed to take any value in $[q,\infty]$, but will be further restricted below. We also refer to the appendix for further notation and results, in particular in the context of tensor product spaces, that will be used in this section.

As described in Section \ref{sec:main_approach} above, the main idea is to synthesize an image via the convolution of a small number of atoms with corresponding coefficient images, where we think of the latter as a sum delta peaks that define the locations where atoms are placed. For this reason, and also due to compactness properties, the coefficient images are modeled as Radon measures in the space $\M(\Omega_\Sigma)$, the dual of $C_0(\Omega_\Sigma)$, where we denote 
\[ \Omega_\Sigma: = \{ x \in \R^d \st \text{there exists }y \in \Sigma \text{ s.t. } x-y \in \Omega \}, \] 
i.e., the extension of $\Omega$ by $\Sigma$. The motivation for using this extension of $\Omega$ is to allow atoms also to be placed arbitrarily close to be boundary (see Figure \ref{fig:lifting_visualization}). We will further use the notation $r'=r/(r-1)$ for an exponent $r \in (1,\infty)$ and denote duality pairings between $L^r$ and $L^{r'}$ and between $\M(\Omega)$ and $C_0(\Omega)$ by $(\cdot,\cdot)$, while other duality pairings (e.g. between tensor spaces) are denoted by $\langle \cdot ,\cdot \rangle$. By $\|\cdot \|_{r}, \|\cdot \|_\M$ we denote standard $L^r$ and Radon norms whenever the domain of definition is clear from the context, otherwise we write $\|\cdot \|_{L^r(\Omega_\Sigma)},\|\cdot \|_{\M(\Omega_\Sigma)} $ etc.

\subsection{The convolutional lasso prior} As first step, we deal with the convolution operator that synthesizes an image from a pair of a coefficient image and an image atom in function space. Formally, we aim to define $K:\M(\Omega_\Sigma) \times L^s(\Sigma) \rightarrow L^q(\Omega)$ as
\[ K(c,p)(x): = \int_{\Omega_\Sigma} p(x-y) \wrt c(y), \]
where we extend $p$ by zero outside of $\Sigma$. An issue with this definition is that, in general, $p$ is only defined Lebesgue almost everywhere and so we have to give a rigorous meaning to the integration of $p$ with respect to an arbitrary Radon measure. To this aim, we define the convolution operator via duality (see \cite{Rudinfourier}). For $c \in \M(\Omega_\Sigma)$, $p \in L^s(\Sigma)$ we define by $K_{c,p}$ the functional on $C(\overline{\Omega})$ as dense subset of $L^{q'}(\Omega)$ as
\[ K_{c,p}(h):= \int _{\R^d} \int _{\R^d} \tilde{h}(z+y) \tilde{p}(z) \wrt z \wrt\tilde{c}(y), \]
where $\tilde{g}$ always denotes the zero extension of the function or measure $g$ outside their domain of definition. Now we can estimate with $D>0$
\[ \left |K_{c,p}(h) \right| \leq  \int _{\R^d} \int _{\R^d} |\tilde{h}(z+y)| |\tilde{p}(z)| \wrt z \wrt|\tilde{c}|(y) \leq \|\tilde{h}\|_{L^{q'}(\R^d)} \|\tilde{p}\|_{L^{q}(\R^d)} \|\tilde{c}\|_{\M(\R^d)} \leq D \|p\|_{s} \|c\|_{\M} \|h\|_{q'}. \]
Hence, by density we can uniquely extend $K_{c,p}$ to a functional in $  L^{q'}(\Omega)^*\simeq L^q(\Omega)$ and we denote by $[K_{c,p}]$ the associated function in $L^q(\Omega)$. Now in case $p$ is integrable w.r.t. $c$ and $x \mapsto \int_{\Omega_\Sigma} p(x-y) \wrt c(y) \in L^q(\Omega)$, we get by a change of variables and Fubini's theorem that for any $h\in C(\overline{\Omega})$
\[ K_{c,p}(h) = \int _{\R^d} \int _{\R^d} \tilde{h}(x) \tilde{p}(x-y) \wrt x \wrt\tilde{c}(y) = \int_\Omega h(x)\left( \int_{\Omega_\Sigma} p(x-y) \wrt c(y) \right) \wrt x.
\]
Hence we get that in this case, $[K_{c,p}](x) = \int_{\Omega_\Sigma} p(x-y) \wrt c(y)$ and defining $K:\M(\Omega_\Sigma) \times L^s(\Omega) \rightarrow L^q(\Omega)$ as
\[ K(c,p):= [K_{c,p}] \]
we get that $K(c,p)$ coincides with the convolution of $c$ and $p$ whenever the latter is well defined. Note that $K$ is bilinear and, as the previous estimate shows, there exists $D>0$ such that $\|K(c,p)\|_q \leq D \|c\|_\M \|p\|_s$. Hence $K \in \B(\M(\Omega_\Sigma) \times L^s(\Sigma),L^q(\Omega))$, the space of bounded bilinear operators (see the appendix).

Using the bilinear operator $K$ and denoting by $k\in \N$ a fixed number for kernels, we now define the convolutional lasso prior for an exponent $s \in [q,\infty]$ and for $u \in L^q(\Omega)$ as

\begin{equation}
\begin{aligned}
 \ncl(u) = 
\inf _ { \substack{  (c_i)_{i=1}^k \subset \M(\Omega_\Sigma) \\   (p_i)_{i=1}^k \subset L^s (\Sigma)}} \, 
  \sum_{i=1}^k \|c_i\|_\M 
 \quad \text{s.t. } 
\left\{
\begin{aligned}
 & \|p_i\|_s \leq 1,\, Mp_i = 0 \quad i=1,\ldots,k,  \\
& u = \sum _{i=1}^k K(c_i,p_i) \text{ in } \Omega,
\end{aligned} \right.
\end{aligned}
\end{equation}
and set $\ncl(u) = \infty$ if the constraint set above is empty. Here, we include an operator $M \in \Lin(L^s(\Sigma),\R^m)$ in our model that optionally allows to enforce additional constraints on the atoms. A simple example of $M$ that we have in mind is an averaging operator, i.e., $Mp:= |\Sigma|^{-1} \int _\Sigma p(x) \wrt x$, hence the constraint that $Mp = 0$ corresponds to a zero-mean constraint. 

\subsection{A convex relaxation}
Our goal is now to obtain a convex relaxation of the convolutional lasso prior. To this aim, we introduce by $\hat{K}$ and $\hat{M}:=I \otimes M$ the lifting of the bilinear operator $K$ and the linear operators $I$ and $M$, with $I\in \Lin(\M(\Omega_\Sigma),\M(\Omega_\Sigma))$ being the identity, to the projective tensor product space $X_s:=\M(\Omega_\Sigma)\ctp L^s(\Sigma) $ (see the appendix). In this space, we consider a reformulation as
\begin{equation} \label{eq:cl_tensor_reformulation}
\begin{aligned}
 \ncl(u) = 
 \inf _ {C \in X_s} \, 
  \|C\|_{ \pi,k,M}\quad  \text{s.t. } 
%\left\{
\begin{aligned}
% & \hat{M}C = 0  \\
& u = \hat{K}C \text{ in } \Omega,
\end{aligned} %\right.
\end{aligned}
\end{equation}
where 
\[ \|C\|_{ \pi,k,M} := \inf \left\{ \sum_{i=1}^k \|c_i\|_\M \|p_i\|_s  \st C = \sum_{i=1}^k c_i \otimes p_i \text{ with } Mp_i = 0 \text{ for } i=1,\ldots,k \right \}. \]
Note that this reformulation is indeed equivalent. Next we aim to derive the convex relaxation of $\ncl$ in this tensor product space. First we consider a relaxation of the functional $\|\cdot \|_{ \pi,k,M}$. To this aim, we need an additional assumption on the constraint set $\ker(M)$, which is satisfied for instance if $s=2$ or for $M=0$, in particular will be fulfilled by the concrete setting we use later on.

\begin{lem}\label{lem:projective_norm_relaxation} Assume that there exists a continuous, linear, norm-one projection onto $\ker(M)$. Then, the convex, lower semi-continuous relaxation of $\|\cdot \|_{\pi,k,M}:X_s \rightarrow \overline{\R}$ is given as
\[ C \mapsto \|C\|_{\pi} + \I_{\ker(\hat{M})}(C), \]
where $\I_{\ker(\hat{M})}(C)= 0$ if $\hat{M}C = 0$ and $\I_{\ker(\hat{M})}(C)= \infty$ else, and $\|\cdot \|_{\pi}$ is the projective norm on $X_s$ given as
\[ \|C\|_\pi = \inf \left \{ \sum_{i=1}^\infty \|c_i\|_\M\|p_i\|_s \st C = \sum_{i=1}^\infty c_i \otimes p_i \right \} .\]
\begin{proof}
At first note that, for a general function $g:X_s \rightarrow \overline{\R}$, its convex, lower-semi continuous relaxation can be computed as the biconjugate $g^{**}:X_s \rightarrow \overline{\R}$, where $g^*(x^*) = \sup_{x\in X_s} \langle x^*,x\rangle  - g(x)$ and $g^{**}(x) = \sup_{x^*\in X^*_s} \langle x^*,x\rangle  - g^*(x^*)$. Keeping this in mind, we note that
\[  \|\cdot \|_{\pi} + \I_{\ker(\hat{M})} \leq  \|\cdot \|_{\pi,k,M}\leq  \|\cdot \|_{\pi,1,M}  , \]
and consequently
\[ \|\cdot \|_{\pi} + \I_{\ker(\hat{M})} \leq  \left( \|\cdot \|_{\pi,k,M} \right)^{**}  \leq  \left( \|\cdot \|_{\pi,1,M}  \right)^{**}.\]
Hence the assertion follows if we show that $\left( \|\cdot \|_{\pi,1,M}  \right)^{**} \leq \|\cdot \|_{\pi} + \I_{\ker(\hat{M})}.$
To this aim, we first show that $\|\cdot \|_{\pi,M} \leq \|\cdot \|_{\pi} + \I_{\ker(\hat{M})} $, where we set $\|C\|_{ \pi,M} = \|C\|_{ \pi,\infty,M}$.
Let $C \in X_s$ be such that $\hat{M}C=0$ and take $(c_i)_i$, $(p_i)_i$ be such that $\|C\|_\pi \geq \sum_{i=1}^\infty \|c_i\|_\M \|p_i\|_s - \epsilon$ for some $\epsilon>0$ and $C = \sum_{i=1}^\infty c_i \otimes p_i$. Then, with $P $ the projection to $\ker(M)$ as in the assumption, we get that
\[ 0 = \hat{M}C =  \sum_{i=1}^\infty c_i \otimes M p_i = \sum_{i=1}^\infty c_i \otimes M (p_i - Pp_i). \]
Now remember that, according to \cite[Theorem 2.9]{Ryan}, we have $(\M(\Omega_\Sigma) \ctp \R^m)^* = \mathcal{B}(\M(\Omega_\Sigma) \times \R^m)$ with the norm $\|B\|_{\mathcal{B}}:= \sup \{ |B(x,y)| \st \|x\|_\M \leq 1, \, \|y\|_2 \leq 1 \}$. Taking arbitrary $\psi \in \M(\Omega)^*$, $\phi \in (\R^m)^*$, we get that $B:(c,p) \mapsto \psi(c)\phi(p) \in \mathcal{B}(\M(\Omega_\Sigma) \times \R^m)$ hence
\begin{align*}
 0 = \hat{B}(\hat{M}C) =  \sum_{i=1}^\infty \hat{B}(c_i \otimes M (p_i - Pp_i))
 & =  \sum_{i=1}^\infty \psi(c_i) \phi( M (p_i - Pp_i)) \\
 & = \phi \left( \sum_{i=1}^\infty \psi(c_i) M (p_i - Pp_i) \right) \\
& = \phi \left( M \left(  \sum_{i=1}^\infty \psi(c_i)  (p_i - Pp_i) \right)\right)
 \end{align*}
 and since $\phi$ was arbitrary it follows that $M \left(  \sum_{i=1}^\infty \psi(c_i)  (p_i - Pp_i)  \right) = 0$. Finally, by closedness of $\range(I-P)$ we get that $\sum_{i=1}^\infty \psi(c_i)  (p_i - Pp_i) = 0$ and, since $\M(\Omega_\Sigma)$ has the approximation property (see \cite[Section VIII.3]{DiestelUhl}), from \cite[Proposition 4.6]{Ryan}, it follows that $\sum_{i=1}^\infty c_i \otimes (p_i - Pp_i) = 0$, hence $C = \sum_{i=1}^\infty c_i \otimes  Pp_i$
and by assumption
$\sum_{i=1}^\infty \|c_i\|_\M \|p_i\|_s \geq \sum_{i=1}^\infty \|c_i\|_\M \|Pp_i\|_s $. Consequently, $\|\cdot \|_{\pi} + \I_{\ker(\hat{M})} \geq \|\cdot \|_{\pi,M} - \epsilon$ and since $\epsilon$ was arbitrary, the claimed inequality follows.

Now we show that $\left( \|\cdot \|_{\pi,M}  \right)^{*} \leq \left(  \|\cdot \|_{\pi,1,M}  \right)^{*}$, from which the claimed assertion follows by the previous estimate and taking the convex conjugate on both sides. To this aim, take $(C_n)_n \subset X_s$ such that
\[\left( \|\cdot \|_{\pi,M} \right)^* (B) = \sup _{C \in X_s} \langle B,C\rangle - \|C\|_{\pi,M} = \lim _n \langle B,C_n\rangle - \|C_n\|_{\pi,M} \]
and take $(c_i^n)_i$, $(p_i^n)_i$ such that $Mp^n_i=0$ for all $n,i$ and 
\[ C_n = \sum _{i=1}^\infty c^n_i \otimes p^n_i \quad \text{and}\quad \sum_{i=1}^\infty \|c_i^n\|_\M \|p_i^n\|_s \leq  \|C_n\|_{\pi,M} + 1/n .\]
We then get
\begin{align*}
\left( \|\cdot \|_{\pi,M} \right)^* (B)
&= \lim _n \langle B,C_n\rangle - \|C_n\|_{\pi,M}  \leq  \lim _n \langle B, \sum _{i=1}^\infty c^n_i \otimes p^n_i\rangle - \sum_{i=1}^\infty \|c_i^n\|_\M \|p_i^n\|_s + 1/n \\
&= \lim _n\lim_m  \langle B, \sum _{i=1}^m c^n_i \otimes p^n_i\rangle - \sum_{i=1}^m \|c_i^n\|_\M \|p_i^n\|_s + 1/n \\
&\leq  \lim _n \sup _m \sup_{ \substack{ (c_i)_{i=1}^m ,(p_i)_{i=1}^m \\ Mp_i = 0} }\langle B, \sum _{i=1}^m c_i \otimes p_i\rangle - \sum_{i=1}^m \|c_i\|_\M \|p_i\|_s + 1/n  \\
&=   \sup _m \sup_{\substack{ (c_i)_{i=1}^m ,(p_i)_{i=1}^m \\ Mp_i = 0}} \sum_{i=1}^m B(c_i,p_i) -  \|c_i\|_\M \|p_i\|_s 
\\
&=  \sup_m  m \sup_{\substack{ c,p \\ Mp = 0}}  \left( B(c,p) -  \|c\|_\M \|p\|_s \right).
\end{align*} 
Now it can be easily seen that the last expression equals $0$ in case $|B(c,p)| \leq \|c\|_\M \|p\|_s$ for all $c,p$ with $Mp = 0 $. In the other case, we can pick $\hat{c},\hat{p}$ with $M\hat{p}=0$ and $r>1$ such that $B(\hat{c},\hat{p}) > r\|\hat{c}\|_\M \|\hat{p}\|_s$ and get for any $\lambda >0$ that
\[ \sup_{\substack{ c,p \\ Mp = 0}}  B(c,p) -  \|c\|_\M \|p\|_s \geq B(\lambda \hat{c},\hat{p})  - \|\lambda\hat{c}\|_\M \|\hat{p}\|_s \geq \lambda (r-1) \rightarrow \infty \text{ as }\lambda \rightarrow \infty .\]
Hence, the last line of the above equation is either 0 or infinity and equals 
\[   \sup_{\substack{ c,p \\ Mp = 0 }} \left( B(c,p) -  \|c\|_\M \|p\|_s \right)  = \sup _C \langle B,C \rangle - \|C\|_{\pi,1,M} =  \left( \|\cdot \|_{\pi,1,M}  \right)^*. \qedhere \]
\end{proof}
\end{lem}
	This result suggests that the convex, lower semi-continuous relaxation of \eqref{eq:cl_tensor_reformulation} will be obtained by replacing $\|\cdot \|_{\pi,k,M}$ with the projective tensor norm $\|\cdot \|_{\pi}$ on $X_s$ and the constraint $\hat{M}C = 0$. Our approach to show this will in particular require us ensure lower semi-continuity of this candidate for the relaxation, which in turn requires us to ensure a compactness property of the sublevelsets of the energy appearing in \eqref{eq:cl_tensor_reformulation} and closedness of the constraints. To this aim, we consider a weak* topology on $X_s$ and rely on a duality result for tensor product spaces (see the appendix), which states that, under some conditions, the projective tensor product $X_s = \M(\Omega_\Sigma) \ctp L^s(\Sigma)$ can be identified with the dual of the so called injective tensor product $C_0(\Omega_\Sigma) \itp L^{s'}(\Sigma)$. Different from one what one would expect from the individual spaces, however, this can only be ensured for the case $s < \infty$ which excludes the space $L^\infty(\Sigma)$ for the image atoms. This restriction will also be required later on in order to show well-posedness of a resulting regularization approach for inverse problems, hence we will henceforth always consider the case that $s \in [q,\infty)$ and use the identification $(C_0(\Omega_\Sigma) \itp L^{s'}(\Sigma))^* \hat{=}\M(\Omega_\Sigma) \ctp L^s(\Sigma)$ (see the appendix).

As first step towards the final relaxation result, and also as crucial ingredient for well-posedness results below, we show weak* continuity of the operator $\hat{K}$ on the space $X_s$.
\begin{lem} Let $s \in [q,\infty)$. Then the operator $\hat{K}:X_s \rightarrow L^q(\Omega)$ is continuous w.r.t. weak* convergence in $ X_s$ and weak convergence in $L^{q'}(\Omega)$. Also, for any $\phi \in C_0(\Omega) \subset L^{q'}(\Omega)$ it follows that $\hat{K}^*\phi \in C_0(\Omega_\Sigma) \itp L^{s'}(\Sigma) $ and, via the identification $C_0(\Omega_\Sigma) \itp L^{s'}(\Sigma) \hat{=} C_0(\Omega_\Sigma,L^{s'}(\Sigma))$ (see the appendix), can be given as $K^*\phi(t)= [ x \mapsto \phi (t + x) ]$.
\begin{proof}
First we note that for any $\psi \in C_c(\Omega)$, the function $\hat{\psi} $ defined as $\hat{\psi}(t)=[ x \mapsto \psi (t+x)]$ (where we extend $\psi$ by 0 to $\R^d$) is contained in $C_c(\Omega_\Sigma,L^{s'}(\Sigma))$. Indeed, continuity follows since by uniform continuity for any $\epsilon > 0$ there exists a $\delta> 0$ such that for any $r\in \R^d $ with $|r| \leq \delta $ and $t \in \Omega_\Sigma$ with $t+r\in \Omega_\Sigma$
\[\| \hat{\psi}(t+r)-\hat{\psi}(t)\|_{s'}=  \left( \int_{\Sigma} |\psi (t +r+  x) - \psi (t  + x)|^{s'}\wrt x \right) ^{1/{s'}} \leq \epsilon |\Sigma|^{1/{s'}} .\]
Also, taking $K \subset \subset \Omega$ to be the support of $\psi$ we get, with $K_\Sigma$ the extension of $K$ by $\Sigma$, for any $t \in \Omega_\Sigma \setminus K_\Sigma$ that $\psi (t + x)=0$ for any $x \in \Sigma $ and hence $\hat{\psi}=0 $ in $L^{s'}(\Sigma)$ and $\hat{\psi}\in C_c(\Omega_\Sigma,L^{s'}(\Sigma))$.

Now for $\phi \in C_0(\Omega_\Sigma)$, taking $(\phi_n)_n \in C_c (\Omega_\Sigma)$ to be a sequence converging to $\phi$, we get that
\[ \|  \hat{\phi} - \hat{\phi_n}\|_{C_0(\Omega_\Sigma,L^{s'}(\Sigma)) } = \sup _t \left( \int_\Sigma | \phi (t+x) - \phi_n (t+x) |^{s'} \wrt x \right)^{1/{s'}}\leq \|\phi - \phi_n \|_\infty |\Sigma|^{1/{s'}} \rightarrow 0 .\]
Thus $\hat{\phi}$ can be approximated by a sequence of compactly supported functions and hence $\hat{\phi} \in C_0(\Omega_\Sigma ,L^{s'}(\Sigma) )$. 
Fixing now $u =c \otimes p\in X_s$ we note that for any $\psi \in C_0(\Omega_\Sigma,L^{s'}(\Sigma))$, the function $t \rightarrow \int_\Sigma \psi(t)(x)p(x)\wrt x$ is continuous, hence we can define the linear functional
 \[ F_{u}(\psi):= \int_{\Omega_\Sigma} \int _\Sigma \psi(t)(x) p(x) \wrt x \wrt c(t) \]
and get that $F_u$ is continuous on $C_0(\Omega_\Sigma,L^{s'}(\Sigma))$. Then, since $\hat{\phi}\in C_0(\Omega_\Sigma) \itp L^{s'}(\Sigma)$ it can be approximated by a sequence of simple tensors $(\sum_{i=1}^{m_n} x_i^n \otimes y_i^n)_n$ in the injective norm, which coincides with the norm in $C_0(\Omega_\Sigma,L^{s'}(\Sigma))$ and, using Lemma \ref{lem:isometriy_itp_cinfty} in the appendix, we get
\begin{align*}
 \langle u,\hat{\phi} \rangle & = \lim_n \langle u,\sum_{i=1}^{m_n} x_i^n \otimes y_i^n \rangle = \lim_n  \sum_{i=1}^{m_n} (c,x_i^n)(p,y_i^n)  = \lim_n \sum_{i=1}^{m_n}  \int_{\Omega_\Sigma} \int _\Sigma x_i^n(a)y_i^n(b)  p(b)\wrt b\wrt c(a) \\
 & = \lim_n F_u ((\sum_{i=1}^{m_n} x_i^n \otimes y_i^n)_n) = F_u (\hat{\phi})  =   \int_{\Omega_\Sigma} \int _\Sigma \phi(a+b) p(b) \wrt b\wrt c(a)  =  (K(c , p),\phi) 
   = (\hat{K}u,\phi) 
 \end{align*}
Now by density of simple tensors in the projective tensor product, it follows that $K^*\phi = \hat{\phi}$. In order to show the continuity assertion, take $(u_n)_n $ weak * converging to some $u \in X_s$. Then by the previous assertion we get for any $\phi \in C_c (\Omega_\Sigma)$ that
\[ (Ku_n,\phi) = \langle u_n,K^* \phi \rangle \rightarrow \langle u,K^*\phi \rangle = (Ku,\phi) ,\]
hence $(Ku_n)_n$ weakly converges to $Ku$ on a dense subset of $L^{s'}(\Omega)$ which, together with boundedness of $(Ku_n)_n$, implies weak convergence.
\end{proof}
 \end{lem}
We will also need weak*-to-weak* continuity of $\hat{M}$, which is shown in the following lemma in a slightly more general situation than needed.
\begin{lem} Take $s \in [q,\infty)$ and assume that $M\in \Lin(L^s(\Sigma), Z)$ with $Z$ a reflexive space and define  $\hat{M}:=I \tp M \in \Lin (X_s ,\M(\Omega_{\Sigma}) \ctp  Z)$, where $I$ is the identity on $\M(\Omega_{\Sigma})$. Then $\hat{M}$ is continuous w.r.t. weak star convergence in both spaces.
\begin{proof}
Take $(u_n)_n \in X$ weak* converging to some $u \in X$ and write $u_n = \lim _k \sum_{i=1}^k x_i^n \otimes y_i^n$. We note that, since $Z$ is reflexive, it satisfies in particular the Radon Nikod\'ym property (see the appendix) and hence $C(\Omega_{\Sigma}) \itp Z^*$ can be regarded as predual of $\M(\Omega_{\Sigma}) \ctp  Z$ and we test with $\phi \otimes \psi \in C(\Omega_{\Sigma}) \itp Z^*$. Then
\begin{align*}
\langle \hat{M}u_n,\phi \otimes \psi \rangle 
&= \lim_k \sum_{i=1}^k ( x_i^n,\phi) ( My_i^n,\psi) = \lim_k \sum_{i=1}^k ( x_i^n,\phi) ( y_i^n,M^*\psi)  \\
&= \langle u_n , \phi \otimes M^* \psi \rangle \rightarrow \langle u , \phi \otimes M^* \psi \rangle = \langle \hat{M}u,\phi \otimes\psi \rangle ,
\end{align*}
where the convergence follows since $M^* \psi \in L^{s'}(\Sigma)$, the predual of $L^s(\Sigma)$, and hence  $\phi \otimes M^* \psi \in C(\Omega_{\Sigma}) \itp L^{q'}(\Sigma) .$
\end{proof}
\end{lem}

Now we can obtain the convex, lower semi-continuous relaxation of $\ncl$.
\begin{lem} With the assumptions of Lemma \ref{lem:projective_norm_relaxation} and $s \in [q,\infty)$, the convex, l.s.c. relaxation of $\ncl$ is given as
\begin{equation} \label{eq:texture_prio_without_nuc}
 \ncl^{**}(u) = \inf_{C \in X_s} \|C\|_\pi \quad \text{s.t. } \begin{cases} u = \hat{K}C \text{ in }\Omega, \\ \hat{M}C = 0. \end{cases} 
 \end{equation}
\proof
Again we first compute the convex conjugate:
\begin{align*}
 \ncl^* (v) 
 &= \sup_u (u,v) - \ncl (u) = \sup _{ C \in X_s } (\hat{K}C,v) - \|C\|_{\pi,k,M} \\
 & = \sup _{  C \in X_s } \langle C,\hat{K}^*v\rangle -   \|C\|_{\pi,k,M} = \|\hat{K}^ * v \|_{\pi,k,M}^*.
\end{align*}
Similarly, we see that $N^*(v) = \left( \| \cdot  \|_\pi + \I_{\ker(\hat{M})} \right)^* (\hat{K}^* v)$, where
\[ N(u) = \inf_{C \in X_s} \|C\|_\pi \quad \text{s.t. } \begin{cases} u = \hat{K}C \text{ in }\Omega \\ \hat{M}C = 0 \end{cases} \]
Now in the proof of Lemma \ref{lem:projective_norm_relaxation}, we have in particularly shown that $ \left( \| \cdot  \|_\pi + \I_{\ker(\hat{M})} \right)^* = \|\cdot \|_{\pi,k,M}^*$, hence, if we show that $N$ is convex and lower semi-continuous, the assertion follows from $N(u) = N^{**}(u) = \ncl^{**}(u)$. To this aim, take a sequence $(u_n)_n $ in $L^q(\Omega)$  converging weakly to some $u$ for which, without loss of generality, we assume that
\[ \lim_n N(u_n) =  \liminf _n N(u_n) < \infty \]
Now with $(C_n)_n$ such that $\|C_n\|_\pi  \leq N(u_n) +n^{-1}$, $\hat{M}C_n = 0$ and $u_n = \hat{K}C_n$ we get that $(\|C_n\|_\pi)_n$ is bounded. Since $X_s$ admits a separable predual (see the appendix), this implies that $(C_n)_n$ admits a subsequence $(C_{n_i})_i$ weak* converging to some $C$. By weak* continuity of $\hat{K}$ and $\hat{M}$ we get that $u = \hat{K}C$ and $\hat{M}C=0$, respectively, and by weak* lower semi-continuity of $\|\cdot \|_\pi$ it follows that
\[ N(u)  \leq \|C\|_\pi \leq \liminf_i \|C_{n_i}\|_\pi \leq \liminf_i N(u_{n_i}) + {n_i}^{-1} \leq \lim_i N(u_{n_i}) = \liminf N(u_n), \]
which concludes the proof. \qedhere
\end{lem}
This relaxation results suggest to use $N(\cdot)$ as in Equation \eqref{eq:texture_prio_without_nuc} as convex texture prior in the continuous setting. 
There is, however, an issue with that, namely that such a functional cannot be expected to penalize the number of used atoms at all. Indeed, taking some $C= \sum_{i=1}^l c_i \otimes p_i$ and assume that $\|C\|_\pi = \sum_{i=1}^l \|c_i\|_\M \|p_i\|_s$. Now note that we can split any summand $c_{i_0} \otimes p_{i_0}$ as follows: Write  $c_{i_0} =c^1_{i_0} + c^2_{i_0} $ with disjoint support such that $\|c_{i_0} \|_\M = \|c^1_{i_0} \|_\M + \|c^2_{i_0}\|_\M$. Then we can re-write 
\[ c_{i_0} \otimes p_{i_0} = c^1_{i_0}\otimes p_{i_0} + c^2_{i_0}\otimes p_{i_0}\]
which gives a different representation of $C$ by increasing the number of atoms without changing the cost of the projective norm. 
Hence, in order to maintain the original motivation of the approach to enforce a limited number of atoms, we need to add an additional penalty on $C$ for the lifted texture prior.

\subsection{Adding a rank penalization}
Considering the discrete setting and the representation of the tensor $C$ as a matrix, the number of used atoms corresponds to the rank of the matrix, for which it is well known that the nuclear norm constitutes a convex relaxation \cite{fazel2002thesis}. This construction can in principle also be transferred to general tensor products of Banach spaces via the identification (see Proposition \ref{prop:nuclear_norm_operator} in the appendix)
\[ C = \sum_{i=1}^\infty x_i \otimes y_i \in X^* \ctp Y^* \quad \longleftrightarrow \quad  T_C \in \Lin(X,Y^*)  \quad \text{where} \quad T_C(x) = \sum_{i=1}^\infty x_i(x)y_i \]
and the norm
\[ \|C\|_\nuc = \|T_C\|_\nuc = \inf \left\{ \sum_{i=1}^\infty \sigma_i \st T_C (x) = \sum_{i=1}^\infty \sigma_i x_i(x)y_i \text{ s.t. } \|x_i\|_{X^*}\leq 1, \|y_i\|_{Y^*} \leq 1 \right \} .\]
It is important to realize, however, that the nuclear norm of operators depends on the underlying spaces and in fact coincides with the projective norm in the tensor product space (see Proposition \ref{prop:nuclear_norm_operator}). Hence adding the nuclear norm in $X_s$ does not change anything and, more generally, whenever one of the underlying spaces is equipped with and $L^1$-type norm we cannot expect a rank-penalizing effect (consider the example of the previous section).

On the other hand, going back to  the nuclear norm of a matrix in the discrete setting, we see that it relies on orthogonality and an inner product structure and that the underlying norm is the Euclidean inner product norm. Hence an appropriate generalization of a rank penalizing nuclear norm needs to be built on a Hilbert space setting. Indeed, it is easy to see that any operator between Banach spaces with a finite nuclear norm is compact, and in particular for any $T \in \Lin(H_1,H_2)$ with finite nuclear norm and $H_1,H_1$ Hilbert spaces, there are orthonormal systems $(x_i)_i$, $(y_i)_i$ and uniquely defined singular values $(\sigma_i)_i$ such that
\[ Tx = \sum_{i=1}^\infty \sigma_i (x,x_i)y_i \quad \text{and in addition}\quad  \|T\|_\nuc = \sum_{i=1}^\infty \sigma _i .\]

Motivated by this, we aim to define an $L^2$-based nuclear norm as extended real valued function on $X_s$ as convex surrogate of a rank penalization. To this aim, we consider from now on the case $s=2$. Remember that the tensor product $X \otimes Y$ of two spaces $X,Y$ is defined as the vector space spanned by linear mappings $x \otimes y$ on the space of bilinear forms on $X \times Y$, which are given as $x \otimes y (B) = B(x,y)$. Now since $L^2(\Omega_\Sigma)$ can be regarded as subspace of $\M(\Omega_\Sigma)$, also $L^2(\Omega_\Sigma) \otimes L^2(\Sigma) $ can be regarded as subspace of $\M(\Omega_\Sigma) \otimes L^2(\Sigma).$ Further, defining for $C \in L^2(\Omega_\Sigma) \otimes L^2(\Sigma)$,
\[ \|C\|_{\pi,L^2\otimes L^2}:= \inf \left \{ \sum_{i=1}^n  \|x_i\|_{2} \|y_i\|_{2} \st C = \sum_{i=1}^n x_i \otimes y_i,\, x_i \in L^2(\Omega_\Sigma), \, y_i \in L^2(\Sigma), \, n \in \N \right\}, \]
we get that $\|\cdot \|_\pi \leq B \|\cdot \|_{\pi,L^2\otimes L^2}$ for a constant $B>0$, hence also the completion $L^2(\Omega_\Sigma) \ctp L^2(\Sigma) $ can be regarded as subspace of $\M(\Omega_\Sigma) \ctp L^2(\Sigma)$. Further, $L^2(\Omega_\Sigma) \ctp L^2(\Sigma) $ can be identified with the space of nuclear operators $\mathcal{N}(L^2(\Omega_\Sigma),L^2(\Sigma))$ as above such that
\[ \|C\|_{\pi,L^2\otimes L^2} = \sum_{i=1}^\infty \sigma_i(T_C) \qquad \text{with } (\sigma_i(T_C))_i \text{ the singular values of }T_C .\]
Using this, and introducing a potential function $\phi:[0,\infty) \rightarrow [0,\infty)$, we define for $C \in X_2$,
\begin{equation} \label{eq:two_norm_extension}
 \|C\|_{\nuc,\phi} := \begin{cases}
 \sum_{i=1}^\infty \phi(\sigma_i(T_C)) &\text{if } C \in L^2(\Omega_\Sigma) \ctp L^2(\Sigma) ,\\
 \infty & \text{ else.}
 \end{cases}
 \end{equation}
 We will mostly focus on the case $\phi(x)=x$, in which $ \| \cdot \|_{\nuc,\phi}$ coincides with an extension of the nuclear norm and can be interpreted as convex relaxation of the rank. However, since we observed a significant improvement in some cases in practice by choosing $\phi$ to be a semi-convex potential function, i.e., a function such that $\phi + \tau |\cdot |^2$ is convex for $\tau$ sufficiently small, we include the more general situation in the theory.

 \begin{rem}[Sparsity and low-rank]\label{rem:nuclear_norm_hilbert}
It is important to note that  $\|C\|_{\nuc,\phi}< \infty $ restricts $C$ to be contained in the smoother space $L^2(\Omega_\Sigma) \ctp L^2(\Sigma)$ and in particular does not allow for simple tensors $ \sum_{i=1}^k c_i \otimes p_i$ with the $c_i$'s being composed of delta peaks. Thus we observe some inconsistency of a rank penalization via the nuclear norm and a pointwise sparsity penalty, which is only visible in the continuous setting trough the regularity of functions. Nevertheless, such an inconsistency has already been observed in the finite dimensional setting in the context of compressed sensing for low-rank AND sparse matrices, manifested via a poor performance of the sum of a nuclear norm and $\ell^1$ norm for exact recovery (see \cite{Oymak15low_rank_sparse}). As a result, there exists many literature on improved, convex priors for the recovery of low-rank and sparse matrices, see for instance \cite{Richard14_sparse_matrix,Chandrasekaran12convex_geometry,Richard13_low_rank_sparse_lifting}. While such improved priors can be expected to be highly beneficial for our setting, the question does not seem to be solved in such a way that that can be readily applied in our setting.

One direct way to circumvent this inconsistency would be to include an additional smoothing operator for $C$ as follows: Take $S \in \Lin(\M(\Omega_\Sigma),\M(\Omega_\Sigma))$ such that $\text{range}(S) \subset L^2(\Omega_\Sigma)$ to be a weak* to weak* continuous linear operator
 and define the operator $\hat{S}:X_2 \rightarrow X_2$ as $\hat{S}:=S \otimes I_{L^2}$, were $I_{L^2}$ denotes the identity in $L^2(\Sigma)$. Then one could alternatively also use $\|SC\|_{\nuc}$ as alternative for penalizing the rank of $C$ while still allowing $C$ to be a general measure. Indeed, in the discrete setting, by choosing $S$ also to be injective, we even obtain the equality $\rank(SC) = \rank(C)$ (where we interpret $C$ and $SC$ as matrices). In practice, however, we did not observe an improvement by including such a smoothing and thus do not include $\hat{S}$ in our model. 
\end{rem}

\subsection{Well-posedness and a cartoon-texture model}
 
Including $\|C\|_{\nuc,\phi}$ for $C \in X_2$ as additional penalty in our model, we ultimately arrive at the following variational texture prior in the tensor product space $X_2:=\M(\Omega_\Sigma)\ctp L^2(\Sigma) $, which is convex whenever $\phi$ is convex, in particular for $\phi(x) = |x|$.
\begin{equation} \label{eq:texture_prior}
 N_\nu(v) = 
 \inf _ {C \in X} \, 
  \nu \|C\|_{ \pi } + (1-\nu)\|C\|_{\nuc,\phi} 
 \qquad \text{s.t. } 
\left\{
\begin{aligned}
 & \hat{M}C = 0 , \\
& v = \hat{K}C \text{ in } \Omega,
\end{aligned} \right.
\end{equation}
where $\nu \in (0,1)$ is a parameter balancing the sparsity and the rank penalty.

In order to employ $N_\nu$ as a regularization term in an inverse problems setting, we need to obtain some lower-semi continuity and coercivity properties. 
As first step, the following lemma, which is partially inspired by techniques used in \cite[Lemma 3.2]{Bredies09_separable}, shows that, under some weak conditions on $\phi$, $ \| \cdot \|_{\nuc,\phi}$ defines a weak* lower semi-continuous function on $X_2$.
\begin{lem} \label{lem:nuc_phi_lower_semi_cont}Assume that $\phi:[0,\infty) \rightarrow [0,\infty)$ is lower semi-continuous, non-decreasing, that 
\begin{itemize}
\item $\phi(x) \rightarrow \infty $ for $x \rightarrow \infty$ and that
\item there exist $\epsilon,\eta>0$ such that $\phi(x) \geq \eta x$ for $0 \leq x < \epsilon$.
\end{itemize} Then the functional  $\|\cdot \|_{\nuc,\phi}:X_2 \rightarrow \overline{R}$ defined as in \eqref{eq:two_norm_extension} is lower semi-continuous w.r.t. weak* convergence in $X_2$.
\begin{proof}
Take $(C_n)_n \subset X_2$ weak* converging to some $C \in X_2$ for which, w.l.o.g., we assume that
\[ \liminf _n \|C_n\|_{\nuc,\phi} = \lim_n \|C_n\|_{\nuc,\phi} .\]
We only need to consider the case that $(\|C_n\|_{\nuc,\phi})_n$ is bounded, otherwise the assertion follows trivially. Hence we can write $T_{C_n}(x) = \sum_{i=1}^\infty \sigma_i^n (x_i^n,x)y_i^n$ such that $\|C_n\|_{\nuc,\phi} = \sum_{i=1}^\infty \phi(\sigma_i^n) $.
Now we aim to bound $(\|C_n\|_{\pi,L^2 \otimes L^2})_n$ in terms of $(\|C_n\|_{\nuc,\phi})_n$. To this aim, first note that the assumptions in $\phi$ imply that for any $\epsilon' >0$ there is $\eta'>0$ such that $\phi(x) \geq \eta' x$ for all $x < \epsilon'$. Also, $\phi(\sigma_i^n) \leq \|C_n\|_{\nuc,\phi}$ for any $i,n$ and via a direct contradiction argument it follows that there exists $\hat{\epsilon}>0$ such that $\sigma_i^n < \hat{\epsilon}$ for all $i,n$. Picking $\hat{\eta}$ such that $\phi(x) \geq \hat{\eta}x$ for all $x < \hat{\epsilon} $ we obtain
\[ \|C_n\|_{\nuc,\phi} = \sum_{i=1}^\infty \phi(\sigma_i^n ) \geq \hat{\eta} \sum_{i=1}^\infty \sigma_i^n = \hat{\eta} \|C_n\|_{\pi,L^2 \otimes L^2} ,\]
hence $(C_n)_n$ is also bounded as sequence in $L^2(\Omega_\Sigma) \ctp L^2(\Sigma)$ and admits a (non-relabeled) subsequence weak* converging to some $\hat{C} \in L^2(\Omega_\Sigma) \ctp L^2(\Sigma)$, with $L^2(\Omega_\Sigma) \itp L^2(\Sigma)$ being the predual space. By the inclusion $C_0(\Omega_\Sigma) \itp L^2(\Sigma) \subset L^2(\Omega_\Sigma) \itp L^2(\Sigma) $ and uniqueness of the weak* limit, we finally get $\hat{C}=C \in L^2(\Omega_\Sigma) \itp L^2(\Sigma) $ and can write $T_Cx = \sum_{i=1}^\infty \sigma_i (x_i,x) y_i$ and $\|C\|_{\nuc,\phi} = \sum_{i=1}^\infty \sigma_i$. By lower semi-continuity of $\|\cdot \|_{\nuc}$ this would suffice to conclude in the case $\phi(x) = x$. For the more general case, we need to show a point-wise lim-inf property of the singular values. To this aim, note that by the Courant-Fischer min-max principle (see for instance \cite[Problem 37]{Brezis}) for any compact operator $T \in \mathcal{L}(H_1,H_2)$ with $H_1,H_2$ Hilbert spaces and $\lambda_k$ the $k$-th singular value of $T$ sorted in descending order we have
\[ \lambda_k = \sup_{\dim(V)=k} \min _{x \in V, \|x\|=1} \|Tx\|_{H_2} .\]
Now consider $k\in \N$ fixed. For any subspace $V$ with $\dim(V)= k$, the minimum in the equation above is achieved, hence we can denote $x_V$ to be a minimizer and define $F_V(T):= \|Tx_V\|_{H_2}$ such that $\lambda_k = \sup_{\dim(V)=k}F_V(T)$. 
Since weak* convergence of a sequence $(T_n)$ to $T$ in $L^2(\Omega_\Sigma) \ctp L^2(\Sigma)$ implies in particular $T_{n}(x) \rightharpoonup T(x)$ for all $x$, by lower semi-continuity of the norm $\|\cdot \|_{H_2}$ it follows that $F_V$ is lower semi-continuous with respect to weak* convergence. Hence this is also true for the function $T \mapsto \lambda_k(T)$ is being the pointwise supremum of a family of lower semi-continuous functional.
Consequently, for the sequence $(T_{C_n})_n$ it follows that $\sigma_k \leq  
 \liminf_n \sigma^n_k$. Finally, by monotonicity and lower semi-continuity of $\phi$ and Fatou's lemma we conclude
 \[ \|T_{C}\|_{\nuc,\phi} = \sum_{k} \phi(\sigma_k) \leq \sum_{k} \phi (\liminf_n \sigma_k^n)= \sum_{k} \liminf_n  \phi (\sigma_k^n) \leq \liminf_n  \sum_k \phi(\sigma_k^n) \leq \liminf_n \|T_{C_n}\|_{\nuc,\phi} .\]

\end{proof}

\end{lem}

The lemma below now establishes the main properties of $N_\nu$ that in particular allow to employ it as regularization term in an inverse problems setting.

\begin{lem} \label{lem:texture_model_properties} The infimum in the definition of \eqref{eq:texture_prior} is attained and $N_\nu:L^q(\Omega) \rightarrow \overline{\R}$ is convex and lower semi-continuous. Further, any sequence $(v_n)_n$ such that $N_\nu(v_n)$ is bounded admits a subsequence converging weakly in $L^q(\Omega)$.
\begin{proof}
The proof is quite standard, but we provide it for the readers convenience. 
Take $(v_n)_n $ to be a sequence such that $(N_\nu(v_n))_n$ is bounded. Then we can pick a sequence $(C_n)_n$ in $X_2$ such that $\hat{M}C_n = 0$, $v_n = \hat{K}C_n$ and 
\[ \nu \|C_n\|_\pi \leq \nu \|C_n\|_\pi + (1-\nu) \|C_n\|_{\nuc,\phi} \leq  N_\nu(v_n) + n^{-1} \]
This implies that $(C_n)_n$ admits a subsequence $(C_{n_i})_i$ weak* converging to some $C \in X_2$. Now by continuity of $\hat{M}$ and $\hat{K}$ we have that $\hat{M}C = 0$ and that $(v_{n_i})_i = (\hat{K}C_{n_i})_i$ is bounded. Hence also $(v_{n_i})_i$ admits a (non-relabeled) subsequence converging weakly to some $v=\hat{K}C$. 
This already shows the last assertion. In order to show lower semi-continuity, assume that $(v_n)_n$ converges to some $v$ and, without loss of generality, that
\[ \liminf_n N_\nu(v_n) = \lim_n N_\nu(v_n) .\]
Now this is a particular case of the argumentation above, hence we can deduce with $(C_n)_n$ as above that
\begin{align*}
 N_\nu(v) 
 &\leq \nu \|C\|_\pi + (1-\nu) \|C\|_{\nuc,\phi} \leq \liminf_i \nu \|C_{n_i}\|_\pi + (1-\nu) \|C_{n_i}\|_{\nuc,\phi} \\
 &\leq \liminf_i  N_\nu (v_{n_i}) +{n_i}^{-1} = \liminf_n N_\nu(v_n) 
\end{align*}
which implies lower semi-continuity. Finally, specializing even more to the case that $(v_n)_n $ is the constant sequence $(v)_n$, also the claimed existence follows.
\end{proof} 
\end{lem}
 
 In order to model a large class of natural images and to keep the number of atoms needed in the above texture prior low, we combine it with a second part that models cartoon-like images. Doing so, we arrive at the following model

\begin{equation} \label{eq:cart_text_problem} \tag{P}
\min _{u ,v\in L^q(\Omega)} 
   \lambda D(Au,f_0) + s_1(\mu)R(u- v) + s_2(\mu)N_\nu(v) 
   \end{equation} 
   where we assume $R$ to be a functional that models cartoon images, $D(\cdot,f_0):Y \rightarrow \overline{\R}$ is a given data discrepancy, $A \in \Lin(L^q(\Omega),Y)$ a forward model and we define the parameter balancing function
  \begin{equation} \label{eq:parameter_balancing_functions}
   s_1(\mu) = 1 - \min(\mu,0), \quad s_2(\mu) = 1 + \max(\mu,0). 
   \end{equation} 
Now we get the following general existence result.
  \begin{prop} \label{prop:main_existence} Assume that $R:L^q(\Omega) \rightarrow \overline{R}$ is convex, lower semi-continuous and that there exists a finite dimensional subspace $U \subset L^q(\Omega)$ such that for any $u\in L^q(\Omega)$, $v \in U^\perp$, $w \in U$,
  \[ \| v \|_q \leq B R(v),\quad \text{and}\quad R(u+w) = R(u) \] with $B>0$ and $U^\perp$ denoting the complement of $U$ in $L^q(\Omega)$. Further assume that $A \in \Lin(L^q(\Omega),Y)$, $D(\cdot,f_0)$ is convex, lower semi-continuous and coercive on the finite dimensional space $A(U)$ in the sense that for any two sequences $(u_n^1)_n$, $(u_n^2)_n$ such that $(u_n^1)_n \subset U$,  $(u_n^2)_n$ is bounded and $(D(A(u_n^1 + u_n^2),f_0))_n$ is bounded, also $(\|A u^1_n\|_q)_n$ is bounded. Then there exists a solution to \eqref{eq:cart_text_problem}.
  \begin{rem}
  Note that for instance in case $D$ satisfies a triangle inequality, the sequence $(u^2_n)$ in the coercivity assumption is not needed, i.e., can be chosen to be zero.
  \end{rem}
  \begin{proof}
    The proof is rather standard and we provide only a short sketch. Take $((u_n,v_n))_n$ a minimizing sequence for \eqref{eq:cart_text_problem}. From Lemma \ref{lem:texture_model_properties} we get that $(v_n)_n$ admits a (non-relabeled) weakly convergent subsequence. Now we split $u_n = u_n^1 + u_n^2 \in U + U^\perp$ and $v_n = v_n^1 + v_n^2 \in U + U^\perp$ and by assumption get that $\|u^2_n - v^2_n\|_q $ is bounded. But since $(\|v_n\|_q)_n$ is bounded, so is $(\|v^2_n\|_q)_n$ and consequently also $(\|u^2_n\|_q)_n$. Now we split again $u_n^1 = u_n^{1,1} + u_n^{1,2} \in \ker(A) \cap U  + (\ker(A) \cap U)^\perp $ and note that also $(u_n^{1,2} + u_n^2,v_n)$ is a minimizing sequence for \eqref{eq:cart_text_problem}. Hence it remains to show that $(u_n^{1,2})_n$ is bounded in order to get a bounded minimizing sequence. To this aim, we note that $(u_n^{1,2})_n \subset (\ker(A) \cap U)^\perp \cap U$ and that $A$ is injective on this finite dimensional space. Hence $\|u_n^{1,2}\|_q \leq B \|Au_n^{1,2}\|_q$ for some $B>0$ and by the coercivity assumption on the data term we finally get that $(\|u_n^{1,2}\|_q)_n$ is bounded. Hence also $(u_n^{1,2} + u_n^2)_n$ admits a weakly convergent subsequence in $L^q(\Omega)$ and by continuity of $A$ as well as lower semi-continuity of all involved functionals existence of a solution follows.
    \end{proof}
  \end{prop}
  
\begin{rem}[Choice of regularization]\label{rem:tgv_regularization} A particular choice of regularization for $R$ in \eqref{eq:cart_text_problem} that we consider in this paper is $R = \TGVat$, with $\TGVat$ the second order total generalized variation functional \cite{bredies2010tgv}, $q \leq d/(d-1)$ and 
 \[Mp:= ( \int_\Sigma p(x) \wrt x,\int_\Sigma p(x_1,x_2) x_1 \wrt x,\int_\Sigma p(x_1,x_2) x_2 \wrt x).\]  Since in this case \cite{Bredies11_inverse,holler14inversetgv}
 \[ \|u\|_q \leq B\TGVat(u) \] with $B>0$ and for all $u \in \mathcal{P}_1(\Omega)^\perp$, the complement of the first order polynomials, and $\TGVat$ is invariant on first order polynomials, the result of Proposition \ref{prop:main_existence} applies.
 \end{rem}
\begin{rem}[Norm-type data terms]We also note that the result of Proposition \ref{prop:main_existence} in particular applies to $D(w,f_0):= \frac{1}{r}\|w-f_0\|_r ^r$ for any $r \in [1,\infty)$ or $D(w,f_0):= \|w-f\|_\infty$, where we extend the norms by infinity to $L^q(\Omega)$ whenever necessary. Indeed, lower semi-continuity of these norms is immediate for both $r\leq q$ and $r>q$ and since the coercivity is only required on a finite dimensional space, it also holds by equivalence of norms.
   \end{rem}

   \begin{rem}[Inpainting] At last we also remark that the assumptions of Proposition \ref{prop:main_existence} also hold for an inpainting data term defined as
   \[ 
   D(w,f_0):= \begin{cases} 
   0 &\text{if } w=f_0 \text{ a.e. on } \omega \subset \Omega \\ 
   \infty & \text{else,}
   \end{cases} \]
   whenever $ \omega$ has nonempty interior. Indeed, lower semi-continuity follows from the fact that $L^q$ convergent sequences admit pointwise convergent subsequences and the coercivity follows from finite dimensionality of $U$ and the fact that $\omega$ has non-empty interior. 
      \end{rem}
      
      \begin{rem}[Regularization in a general setting] We also note that Lemma \ref{lem:texture_model_properties} provides the basis for employing either $N_\nu$ directly or its infimal-convolution with a suitable cartoon prior as in Proposition \ref{prop:main_existence} for the regularization of general (potentially non-linear) inverse problems and with multiple data fidelities, see for instance \cite{Hofmann07_nonlinear_tikhonov_banach,holler17coupled_reg} for general results in that direction.      
      \end{rem}

\section{The model in a discrete setting} \label{sec:discrete_setting}
This section deals with the discretization of the proposed model and its numerical solution. For the sake of brevity, we provide only the main steps and refer to the publicly available source code \cite{convex_lifting_code} for all details.

We define $U = \R^{N \times M}$ to be the space of discrete grayscale images, $W = \R^{(N+n-1) \times (M+n-1)}$ to be the space of coefficient images and $Z =  \R^{n \times n}$ to be the space of image atoms for which we assume $n<\min\{N,M\}$ and, for simplicity, only consider a square domain for the atoms.
The tensor product of a coefficient image $c\in W$ and a atom $p  \in Z$ is given as $(c \otimes p )_{i,j,r,s} = c_{i,j}p_{r,s} $ and the lifted tensor space is given as the four-dimensional space $X =  \R^{(N+n-1) \times (M+n-1) \times n \times n}$.

\noindent \textbf{Texture norm.} The forward operator $K$ being the lifting of the convolution $c \ast p$ and mapping lifted matrices to the vectorized image space is then given as
\[ (KC)_{i,j} = \sum_{r,s=1}^{n,n}  C_{i+n-r,j+n-s,r,s} \]
and we refer Figure \ref{fig:lifting_visualization} for a visualization in the one dimensional case. Note that by extending the first two dimensions of the tensor space to $N+n-1$, $M+n-1$ we allow to place an atom at any position where it still effects the image, also partially outside the image boundary.

Also we note that, in order to reduce dimensionality and accelerate the computation, we introduce a stride parameter $\eta \in \N$ in practice which introduces a stride on the possible atom positions. That is, the lifted tensor space and forward operator is reduced in such a way that the grid of possible atom positions in the image is essentially $\{ (\eta i,\eta j) \st i,j\in \N, (\eta i,\eta j) \in \{1,\ldots,N\} \times \{1,\ldots,M\} \} $. This reduces the dimension of the tensor space by a factor $\eta^{-2}$, while for $\eta>1$ it naturally does not allow for arbitrary atom positions anymore and for $\eta=n$ it corresponds to only allowing non-overlapping atoms. In order to allow for atoms being placed next to each other it is important to choose $\eta$ to be a divisor of the atom-domain size $n$ and we used $n=15$ and $\eta=3$ in all experiments of the paper. In order to avoid extensive indexing and case distinctions, however, we only consider the case $\eta=1$ here and refer to the source code \cite{convex_lifting_code} for the general case.

As a straightforward computation shows that, in the discrete lifted tensor space, the projective norm corresponding to discrete $\|\cdot \|_1$ and $\|\cdot \|_2$ norms for the coefficient images and atoms, respectively, is given as a mixed 1-2 norm as
\[ \|A\|_\pi = \|A\|_{1,2} = \sum_{i,j=1}^{N,M} \sqrt{ \sum_{r,s=1}^{n,n} A_{i,j,r,s}^2} .\]
The nuclear norm for a potential $\phi$ on the other hand reduces to the evaluation of $\phi$ on the singular values of a matrix-reshaping of the lifted tensors and is given as 
\[ \|A\|_{\nuc,\phi} = \sum_{i=1}^{nn} \phi(\sigma_i), \quad \text{with }(\sigma_i)_i\text{ the singular values of }B = \reshape{A}.\]
where $\reshape{A}$ denotes a reshaping of the tensor $A$ to a matrix of dimensions $NM \times nn$. 
For the potential function $\phi$ we consider two choices: Mostly we are interested in $\phi(x) = x$ which yields a convex texture model and enforces sparsity of the singular values. A second choice we consider is $\phi:[0,\infty) \rightarrow [0,\infty)$ given as\\
\begin{minipage}[l]{0.49\linewidth}
\centering
\begin{tikzpicture}
      \draw[->] (0,0) -- (6*0.75,0) node[right] {$x$};
      \draw[->] (0,0) -- (0,1) node[above] {$\phi(x)$};
      \draw[scale=6.0,domain=0:0.25,smooth,variable=\x,blue] plot ({\x},{\x - 2.0*0.99*\x*\x});
      \draw[scale=6.0,domain=0.25:0.75,smooth,variable=\y,blue] plot ({\y},{0.01*\y + 0.99/8});
\end{tikzpicture}
\end{minipage}
\begin{minipage}[r]{0.49\linewidth}
\centering
\begin{equation} \label{eq:semiconv_potential}
 \phi(x) = \begin{cases} 
x - \epsilon \delta x^2 & x \in [0,\frac{1}{2\epsilon}] \\
(1-\delta) x + \frac{\delta}{4 \epsilon} & \text{else,}\end{cases}
\end{equation}
\end{minipage}
\vspace*{0.5\baselineskip}

where $\delta < 1$, $\delta \approx 1$ and $\epsilon >0$. 
It is easy to see that $\phi $ fulfills the assumptions of Lemma \ref{lem:nuc_phi_lower_semi_cont} and that $\phi$ is semi-convex, i.e., $\phi + \rho |\cdot |^2$ is convex for $\rho>\delta\epsilon$. While the results of Section \ref{sec:continuous_setting} hold for this setting even without the semi-convexity assumption, we can in general not expect to obtain an algorithm that provably delivers a globally optimal solution in the semi-convex (or generally non-convex) case. The reason for using a semi-convex potential rather than a arbitrary non-convex one is twofold: First, for a suitably small stepsize $\tau$ the proximal mapping
\[ \prox_{\tau,\phi}(\hat{u}) = \argmin _u \frac{\|u-\hat{u}\|_2 ^2}{2\tau} + \phi(u) \]
is well defined and hence proximal-point type algorithms are applicable at least conceptually. Second, since we employ $\phi$ on the singular values of the lifted matrices $A$, it will be important for numerical feasibility of the algorithm that the corresponding proximal mapping on $A$ can be reduced to a proximal mapping on the singular values. While this is not obvious for a general choice of $\phi$, it is true (see Lemma \ref{lem:semiconvex_prox} blow) for semi-convex $\phi$ with suitable parameter choices.

\noindent\textbf{Cartoon prior} As cartoon-prior we employ the second-order total generalized variation functional which we define for fixed parameters $(\alpha_0,\alpha_1) = (\sqrt{2},1) $ and a discrete image $u \in U$ as
\[ \TGVat(u) = \min_{v \in U^2} \alpha_1 \|\nabla u - v \|_1 + \alpha_0 \|\symgrad v\|_1. \]
Here $\nabla$ and $\symgrad $ denote discretized gradient and symmetrized Jacobian operators, respectively and we refer to \cite{holler15tgvrec_p2} and the source code \cite{convex_lifting_code} for details on a discretization of $\TGVat$. To ensure a certain orthogonality of the cartoon and texture part, we further define the operator $M$ that incorporates atom-constraints, to evaluate the 0th and 1st moments of the atoms, which in the lifted setting yields
\[ (MC)_{i,j}:= \left( \sum_{r,s=1}^{n,n} C_{i,j,r,s},\sum_{r,s=1}^{n,n} r C_{i,j,r,s}, \sum_{r,s=1}^{n,n} s C_{i,j,r,s}\right) .\]

The discrete version of \eqref{eq:cart_text_problem} is then given as
\begin{equation} \label{eq:discrete_problem_general} \tag{DP}
 \min _{  \substack{ u\in U ,C \in X \\ MC = 0 }} \lambda D(Au,f_0)  + s_1(\mu) \TGVat(u-KC) + s_2(\mu)\left(\nu\|C\|_{1,2} + (1-\nu)\|C\|_{\nuc,\phi} \right) , 
 \end{equation}
where the parameter balancing functions $s_1,s_2$ are given as in \eqref{eq:parameter_balancing_functions} and the model depends on three parameters $\lambda,\mu,\nu$, with $\lambda$ defining the trade-off between data and regularization, $\mu$ defining the trade-off between the cartoon and the texture part and $\nu$ defining the trade-off between sparsity and low-rank of the tensor $C$.

\noindent\textbf{Numerical solution.} For the numerical solution of \eqref{eq:discrete_problem_general} we employ the primal-dual algorithm of \cite{pock2011primaldual}. Since the concrete form of the algorithm depends on whether the proximal mapping of the data term $u \mapsto D(Au,f_0)$ is explicit or not, we replace the data term $D(Au,f_0)$ by
\[ D_1(Au,f_0) + D_2(u,f_0) \] where we assume the proximal mappings of $v \mapsto D_i(v,f_0)$ to be explicit and, depending on the concrete applications, set either $D_1$ or $D_2$ to be the constant zero function. 

Denoting by $g^*(v):= \sup_{w} (v,w) -g(w)$ the convex conjugate of a function $g$, with $(\cdot ,\cdot)$ being the standard inner product of the sum of all pointwise-products of entries of $v$ and $w$, we reformulate \eqref{eq:discrete_problem_general} to a saddle-point problem as
\begin{align*}
\eqref{eq:discrete_problem_general}
 & \Leftrightarrow \min _{  u\in U ,C \in X} & \lambda D_1(Au,f_0) + \lambda D_2(u,f_0)  + s_1(\mu) \TGVat(u-KC) + s_2(\mu)(\nu\|C\|_{1,2}  \\ 
 & &  + (1-\nu)\|C\|_{\nuc,\phi} )   + \I_{\ker(M)}(C) \\
 & \Leftrightarrow \min _{  \substack{ u\in U ,C \in X \\ v \in U^ 2 }} \max_{p,q,d,r,m} 
  & (Au,d) - (\lambda D_1(\cdot,f_0))^*(d) +  (\nabla (u - KC) - v,p) - \I_{\|\cdot \|_\infty \leq \alpha_1s_1(\mu) }(p)\\
  &&  + (\symgrad v,q)  - \I_{\|\cdot \|_\infty \leq \alpha_0s_1(\mu) }(q) +  (C,r) - (s_2(\mu)\nu\|\cdot \|_{1,2})^*(r) \\
  &&+ (MC,m) -  \I_{\{0\}}^*(m)  + \lambda D_2(u,f_0) + s_2(\mu) (1-\nu) \|C\|_{\nuc,\phi} \\
  & \Leftrightarrow \min _{x = (u,v,C)} \max_{y = (p,q,d,r,m)} &(Bx,y) - F^*(y) + G(x).
\end{align*}
Here, the dual variables $(p,q,d,r,m) \in (U^2,U^3,A(U),U,U^3)$ are in the image space of the corresponding operators, $\I_S(z) = 0$ if $z \in S$ and $\I_S(z) = \infty $ else, $\{ \|\cdot \|_\infty \leq \delta \}:= \{ z\st \|z\|_\infty \leq \delta \}$ with $\|z  \|_\infty = \|(z_1,\ldots,z_l)  \|_\infty = \sup_{i,j} \sqrt{\sum_{s=1}^l (z_{i,j}^s)^2}$ a point-wise infinity norm on $z \in U^l$. The operator $B$ and the functional $G$ are given as 
\[B(u,v,C) = (\nabla u - \nabla KC - v,\symgrad v,Au,C,MC), \quad G(x) = G(u,v,C) = \lambda D_2(u,f_0) + s_2(\mu) (1-\nu) \|C\|_{\nuc,\phi}\]
 and $F^*(y) = F^*(p,q,d,r,m)$ summarizes all the conjugate functionals as above. Applying the algorithm of \cite{pock2011primaldual} to this reformulation yields the numerical scheme as in Algorithm \ref{alg:algorithm_general}.

\begin{algorithm}[t]
  \begin{algorithmic}[1]
\onehalfspacing

\Function{cart\_text\_recon}{$f_0$}

\State $(u,v,C) \gets  (0,0,0), (p,q,d,r,m) \gets (0,0,0,0,0)$
\State choose $\sigma, \tau > 0$ 
\Repeat
	\State \textbf{Dual updates}
	\State $ p \gets \proj_{\alpha _1} \left( p + \sigma (\nabla (\overline{u} - K \overline{C})-\overline{v}) \right)$
	\State $q \gets \proj_{\alpha _0} (q + \sigma \mathcal{E}\overline{v})$
	\State $d \gets \prox_{\sigma,(\lambda D_1(\cdot,f_0))^*} (d + \sigma A\overline{u})$
	\State $r \gets \prox_{\sigma,(s_2(\mu)\nu\|\cdot \|_{1,2})^*} (r + \sigma \overline{C})$
	\State $m \gets  (m + \sigma M\overline{C})$
		\State \textbf{Primal updates}
	\State $u_+ \gets \prox_{\tau,(\lambda D_2(\cdot,f_0))^*}(u - \tau ( \nabla^* p + A^*d))$
	\State $v_+ \gets v - \tau (-p + \symgrad ^* q)$
	\State $C_+ \gets \prox_{\tau,s_2(\mu)(1-\nu)\|\cdot \|_{\nuc,\phi}} ( C - \tau ( - K^*\nabla ^* p + r + M^*m)$	
		\State \textbf{Extrapolation and update}
	\State $(\overline{u},\overline{v},\overline{C}) \gets 2(u_+,v_+,C_+) - (u,v,C)$, 
	\State $(u,v,C) \gets (u_+,v_+,C_+)$
\Until{Stopping criterion fulfilled}
\State \Return{$(u_+,KC_+,\text{RSV}(C_+)$)}
\EndFunction
\end{algorithmic}
\caption{Scheme of implementation for the numerical solution of \eqref{eq:discrete_problem_general}}\label{alg:algorithm_general}
\end{algorithm}
    
Note that there, we set either $D_1(\cdot,f_0) \equiv 0 $ such that the dual variable $d$ is constant 0 and line $8$ of the algorithm can be skipped, or we set $D_2(\cdot,f_0) \equiv 0$ such that the proximal mapping in line $12$ reduces to the identity. The concrete choice of $D_1,D_2$ and the proximal mappings will be given in the corresponding experimental sections.
 All other proximal mappings can be computed explicitly and reasonably fast: The mappings $\proj_{\alpha _1}$ and $\proj_{\alpha _1}$ can be computed as pointwise projections to the $L^\infty$-ball (see for instance \cite{holler15tgvrec_p2}) and the mapping $\prox_{\sigma,(s_2(\mu)\nu\|\cdot \|_{1,2})^*}$ is a similar projection given as
\[ \prox_{\sigma,(s_2(\mu)\nu\|\cdot \|_{1,2})^*}(R)_{i,j,l,s} = R_{i,j,l,s}/ \left( \max \left\{1,(\sum_{l,s=1}^{n,n} R_{i,j,l,s}^2)^{1/2}/ (s_2(\mu)\nu) \right\} \right) .\]
Most of the computational effort lies in the computation of $\prox_{\tau,s_2(\mu)(1-\nu)\|\cdot \|_{\nuc,\phi}}$, which, as the following lemma shows, can be computed via an SVD and a proximal mapping on the singular values.

\begin{lem} \label{lem:semiconvex_prox} Let $\phi:[0,\infty) \rightarrow [0,\infty)$ be a differentiable and increasing function and $\tau,\rho>0$ be such that $x \mapsto \frac{x^2}{2\tau}+ \rho \phi(x)$ is convex on $[0,\infty)$. Then the proximal mapping of $\rho\|\cdot \|_{\nuc,\phi}$ for parameter $\tau$ is given as
\[ \prox_{\tau,\rho \|\cdot \|_{\nuc,\phi}} (A) = \reshape[(N,M,n,n)]{ \left( U\diag((\prox_{\tau,\rho\phi}(\sigma_i))_i)V^T \right)} \]
where $\reshape{A} = U \Sigma V^T$ is the SVD of $\reshape{A}$  and for $x_0 \geq 0$
\[ \prox_{\tau,\rho\phi}(x_0) = \min_{x} \frac{|x-x_0|^2}{2\tau} + \rho\phi(|x|) .\]
In particular, in case $\phi(x) = x$ we have
\[ \prox_{\tau,\rho \phi}(x_0) =  
\begin{cases}
0 & \text{if }  \quad  0 \leq x_0 \leq \tau \rho, \\
x_0 - \tau\rho  & \text{else,}\\
\end{cases} \]
and in case
\[ \phi(x) = \begin{cases} 
x - \epsilon \delta x^2 &\text{if } x \in [0,\frac{1}{2\epsilon}], \\
(1-\delta) x + \frac{\delta}{4 \epsilon} & \text{else,}
\end{cases}
\]
we have that that $x \mapsto \frac{x^2}{2\tau}+ \rho \phi(x)$ is convex whenever $\tau  \leq \frac{1}{2\epsilon \delta \rho}$ and in this case
\[ \prox_{\tau,\rho \phi}(x_0) = \begin{cases}
 0 & \text{if }  \quad  0 \leq x_0 \leq \tau \rho,\\
\frac{x_0 - \tau\rho}{1 - 2 \epsilon \delta \tau \rho} &\text{if } \quad \tau\rho < x_0  \leq \frac{1}{2 \epsilon} + \tau\rho(1-\delta), \\
x_0 - \tau\rho(1-\delta) &\text{if } \quad   \frac{1}{2 \epsilon} + \tau\rho(1-\delta)  < x_0 .\\
\end{cases} \]
\begin{proof} 
At first note that it suffices to consider $\rho \|\cdot \|_{\nuc,\phi}$ as a function on matrices and show the assertion without the reshaping operation. For any matrix $A$, we denote by $A = U_A \Sigma_A V^T _A$ the SVD of $A$ and $\Sigma_A = \text{diag}((\sigma^A_i)_i)$ contains the singular values sorted in non-increasing order, where $\Sigma_A$ is uniquely determined by $A$ and $U_A,V_A$ are chosen to be suitable orthonormal matrices.

We first show that $G(A):= \frac{\|A\|^2_2}{2 \tau} + \rho \|A \|_{\nuc,\phi}$ is convex. For $\lambda \in [0,1]$, $A,B$ matrices we get by sub-additivity of the singular values (see for instance \cite{Thompson75svd_sum}) that
\begin{align*}
G(\lambda A + (1-\lambda) B)
& =  \sum_{i} \frac{1}{2\tau}(\sigma_i^{\lambda A + (1-\lambda ) B } )^2 + \rho \phi(\sigma_i^{\lambda A + (1-\lambda ) B}) \\
& \leq   \sum_{i} \frac{1}{2\tau}(\lambda \sigma^A_i + (1-\lambda ) \sigma^B_i )^2 + \rho \phi(\lambda \sigma^A_i + (1-\lambda ) \sigma^B_i)\\
& \leq   \sum_{i} \frac{\lambda}{2\tau} (\sigma^A_i )^2 + \frac{1-\lambda}{2\tau}(\sigma^B_i)^2 + \rho\lambda \phi(\sigma^A_i) + \rho(1-\lambda) \phi(\sigma^B_i)\\
& \leq \lambda G(A) + (1-\lambda) G(B).
\end{align*}
Now with $H(A):=\frac{\|A-A_0\|^2_2}{2 \tau} + \rho \|A \|_{\nuc,\phi}$ we get that $H(A) = G(A)  - \frac{1}{2\tau}(2(A,A_0)  + \|A_0\|_2^2)$, thus also $H$ is convex.
Hence first order optimality conditions are necessary and sufficient and we get (using the derivative of the singular values as in \cite{Papadopoulo00_singular_values}) with $DH $ the derivative of $H$ that $A = \prox_{\tau,\rho \|\cdot \|_{\nuc,\phi}} (A_0)$ is equivalent to 
\begin{align*}
 0 & = DH(A) = (A-A_0) + \tau \rho U_A \diag( (\phi^\prime( \sigma^A_i))_i) V_A^T \\
 & = -A_0 + U_A ( \Sigma_A + \tau \rho \diag( (\phi^\prime( \sigma^A_i))_i) )V_A^T
\end{align*}
and consequently to
\[ \sigma^{A_0}_i = \sigma^A_i + \tau \rho \phi^\prime (\sigma^A_i) \]
which is equivalent to 
\[ \sigma ^ {A_i} = \prox _{\tau,\rho\phi} (\sigma_i^{A_0}) \]
as claimed. The other results follow by direct computation.
\end{proof}
\end{lem}

Note also that, in Algorithm \ref{alg:algorithm_general}, $KC_+$ returns the part of the images that is represented by the atoms (the ``texture part'') and $\text{RSV}(C_+)$ stand for right-singular values of $\reshape{(C_+)}$ and returns the image atoms.  For the sake of simplicity, we use a rather high, fixed number of iterations in all experiment but note that, alternatively, a duality-gap based stopping criterion (see for instance \cite{holler15tgvrec_p2}) could be used.

\begin{figure}[t]
\center 

\includegraphics[width=0.24\linewidth]{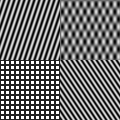}
\includegraphics[width=0.24\linewidth]{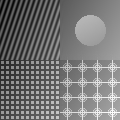}
\includegraphics[width=0.24\linewidth]{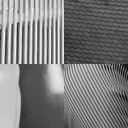}
\includegraphics[width=0.24\linewidth]{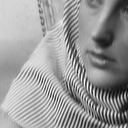}
\caption{\label{fig:test_images} Different test images we will refer to as: Texture, Patches, Mix, Barbara}

\end{figure}

\section{Numerical results} \label{sec:numerical_results}
In this section we present numerical results obtained with the proposed method as well as its variants and compare to existing methods. We will mostly focus on the setting of \eqref{eq:discrete_problem_general}, where $\phi(x) = |x|$ and we use different data terms $D$. 
Hence, the regularization term is convex and consists of $\TGVat$ for the cartoon part and a weighted sum of a nuclear norm and $\ell^{1,2}$ norm for the texture part. Besides this choice of regularization (called \ctc), we will compare to pure $\TGVat$ regularization (called \tgv), the setting of \eqref{eq:discrete_problem_general} with the semi-convex potential $\phi$ as in \eqref{eq:semiconv_potential} (call \ctsc) and the setting of \eqref{eq:discrete_problem_general} with $\TGV$ replaced by $\I_{\{0\}}$, i.e., only the texture norm is used for regularization, and $\phi(x) = |x|$ (called \txt). Further, in the last subsection, we also compare to other methods as specified there. For \ctc{} and \ctsc{} we use the algorithm described in the previous section (where convergence can only be ensured for \ctc{}) and for the other variants we use a adaption of the algorithm to the respective special case.

We fix the size of the atom domain to $15 \times 15$ pixel and the stride to $3$ pixel (see Section \ref{sec:discrete_setting}) for all experiments, and use four different test images (see Figure \ref{fig:test_images}): The first two are synthetic images of size $120\times 120$, containing four different blocks of size $60\times 60$, whose size is a multiple of the chosen atom-domain size. The third and fourth image have size $128 \times 128$ (not being a multiple of the atom-domain size) and the third image contains four sections of real images of size $64 \times 64 $ each (again not a multiple of the atom-domain size). All but the first image contain a mixture of texture and cartoon parts. The first four subsections consider only convex variants of our method ($\phi(x) = |x|$) and the last one considers the improvement obtained with a non-convex potential $\phi$ and also compares to other approaches.

Regarding the choice of parameters for all methods, we generally aimed to reduce the number of varying parameters for each method as much as possible such that for each method and type of experiment, at most two parameters need to be optimized.
Whenever we incorporate the second order TGV functional for the cartoon part, we fix the parameters $(\alpha_0,\alpha_1)$ to $(\sqrt{2},1)$. The method \ctc{} then essentially depends on the three parameters $\lambda, \mu,\nu$. We experienced that the choice of $\nu$ is rather independent of the data and type of experiments, hence we leave it fixed for all experiments with incomplete or corrupted data, leaving our method with two parameters to be adapted: $\lambda$ defining the tradeoff between data and regularization and $\mu$ defining the tradeoff between cartoon and texture regularization. For the semi-convex potential we choose $\nu$ as with the convex one, fix $\delta =0.99$ and use two different choice of $\epsilon$, depending on the type of experiment, hence again leaving two parameters to be adapted. A summary of the parameter choice for all methods is provided in Table \ref{tbl:parameter_choice} below.

We also note that, whenever we tested a range of different parameters for any method presented below, we show the visually best results in the figure. Those are generally not the ones delivering the best result in terms of peak-signal-to-noise ratio, and for the sake of completeness we also provide in Table \ref{tbl:psnr_results} the best PSNR result obtained with each method and each experiment over the range of tested parameters.

\subsection{Image-atom learning and texture separation}
\begin{figure}[t]
\center 

\hspace*{0.3\linewidth}
\includegraphics[width=0.3\linewidth]{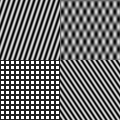}
\includegraphics[width=0.3\linewidth]{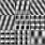}

\includegraphics[width=0.3\linewidth]{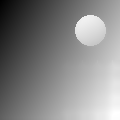}
\includegraphics[width=0.3\linewidth]{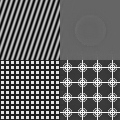}
\includegraphics[width=0.3\linewidth]{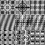}

\includegraphics[width=0.3\linewidth]{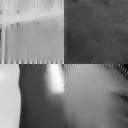}
\includegraphics[width=0.3\linewidth]{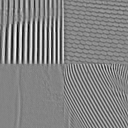}
\includegraphics[width=0.3\linewidth]{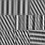}

\includegraphics[width=0.3\linewidth]{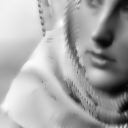}
\includegraphics[width=0.3\linewidth]{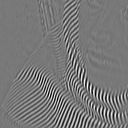}
\includegraphics[width=0.3\linewidth]{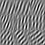}

\caption{\label{fig:cart_text_decomp} Cartoon-texture decomposition (rows 2-4) and nine most important learned atoms for different test images and the methods \txt{} (row 1) and \ctc{} (rows 2-4).}

\end{figure}

As first experiment we test the method \ctc{} for learning image atoms and texture separation directly on the ground truth images. To this aim, we use 
\[ D_1 \equiv 0, \qquad D_2(u,f_0) = \I_{ \{ 0\} } (u-f_0) .\] 
The results can be found in Figure \ref{fig:cart_text_decomp}, where for the pure texture image we used only the texture norm (i.e. the method \txt) without the TGV part for regularization. 

It can be observed that the proposed method achieves a good decomposition of cartoon and texture and also is able to learn the most important image structure effectively. While there are some repetitions of shifted structures in the atoms, the different structures are rather well-separated and the first nine atoms corresponding to the nine largest singular values still contain the most important features of the texture parts.

\subsection{Inpainting and leaning from incomplete data}
This section deals with the task of inpainting a partially available image and learning image atoms from this incomplete data. For reference, we also provide results with pure $\TGVat$ regularization (the method \tgv). The data fidelity in this case is
\[ D_1 \equiv 0, \qquad D_2(u,f_0) = \I _{\{v \st v_{i,j} = (f_0)_{i,j} \text{ for } (i,j) \in \mathcal{M}\}}  (u),\]
with $\mathcal{M}$ the index set of known pixels. Again we use only the texture norm for the first image (the method \txt) and the cartoon-texture functional for the others.

\begin{figure}[t]
\center 

\includegraphics[width=0.24\linewidth]{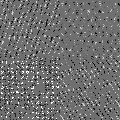}
\includegraphics[width=0.24\linewidth]{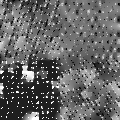}
\includegraphics[width=0.24\linewidth]{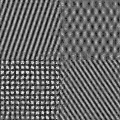}
\includegraphics[width=0.24\linewidth]{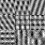}

\includegraphics[width=0.24\linewidth]{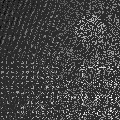}
\includegraphics[width=0.24\linewidth]{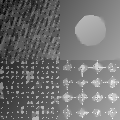}
\includegraphics[width=0.24\linewidth]{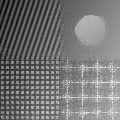}
\includegraphics[width=0.24\linewidth]{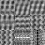}

\includegraphics[width=0.24\linewidth]{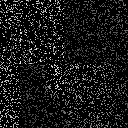}
\includegraphics[width=0.24\linewidth]{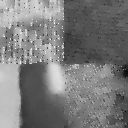}
\includegraphics[width=0.24\linewidth]{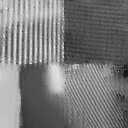}
\includegraphics[width=0.24\linewidth]{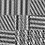}

\includegraphics[width=0.24\linewidth]{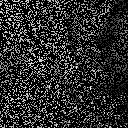}
\includegraphics[width=0.24\linewidth]{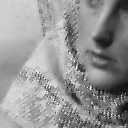}
\includegraphics[width=0.24\linewidth]{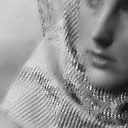}
\includegraphics[width=0.24\linewidth]{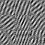}

\caption{\label{fig:inpainting_test} Image inpainting from incomplete data. From left to right: Data, TGV-based reconstruction, proposed method (only \txt{} in first row),  nine most important learned atoms.  Rows 1,3: 20\% of pixels, rows 2,4: 30\% of pixels.}

\end{figure}

The results can be found in Figure \ref{fig:inpainting_test}. For the first and third image, $20\%$ of the pixels where given while for the other two, $30\%$ were given. It can be seen that our method is generally still able to identify the underlying pattern of the texture part and to reconstruct it reasonably well. Also the learned atoms are reasonable and are in accordance with the ones learned from the full data as in the previous section. In contrast to that, pure $\TGV$ regularization (which assumes piecewise smoothness) has no chance to reconstruct the texture patterns. For the cartoon part, both methods are comparable. It can also be observed that the target-like structure in the bottom right of the second image is not reconstructed well and also not well identified with the atoms (only the 8th one contains parts of this structure). The reason might be that due to the size of the repeating structure there is not enough redundant information available to reconstruct it from the missing data. Concerning the optimal PSNR values of Table \ref{tbl:psnr_results}, we can observe a rather strong improvement with \ctc{} compared to \tgv.

\subsection{Learning and separation under noise}

In this section we test our method for image-atom-learning and de-noising with data corrupted by Gaussian noise (with standard deviation 0.5 and 0.1 times the image range for the Texture and the other images, respectively). Again we compare to $\emph{TGV}$ regularization in this section (but also to other methods in Section \ref{sec:comparison} below) and use the texture norm for the first image (the method \txt). The data fidelity in this case is
\[ D_1 \equiv 0, \qquad D_2(u,f_0) = \frac{1}{2}\|u-f_0\|_2^2.\]

It can be observed that also under the presence of rather strong noise, our method is able to learn some of the main features of the image within the learned atoms. Also the quality of the reconstructed image is improved compared to \tgv, in particular for the right-hand side of the \emph{Mix} image, where the top left structure is only visible in the result obtained with \ctc. On the other hand, the circle of the \emph{Patches} image obtained with \ctc{} contains some artifacts of the texture part. Regarding the optimal PSNR values of Table \ref{tbl:psnr_results}, the improvement with \ctc{} compared to \tgv{} is still rather significant.

\begin{figure}[t]
\center 

\includegraphics[width=0.24\linewidth]{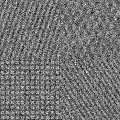}
\includegraphics[width=0.24\linewidth]{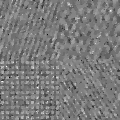}
\includegraphics[width=0.24\linewidth]{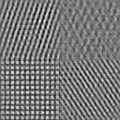}
\includegraphics[width=0.24\linewidth]{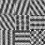}

\includegraphics[width=0.24\linewidth]{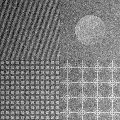}
\includegraphics[width=0.24\linewidth]{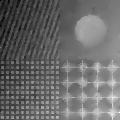}
\includegraphics[width=0.24\linewidth]{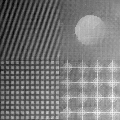}
\includegraphics[width=0.24\linewidth]{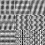}

\includegraphics[width=0.24\linewidth]{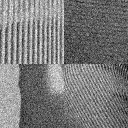}
\includegraphics[width=0.24\linewidth]{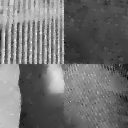}
\includegraphics[width=0.24\linewidth]{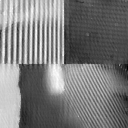}
\includegraphics[width=0.24\linewidth]{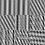}

\includegraphics[width=0.24\linewidth]{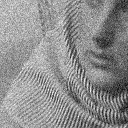}
\includegraphics[width=0.24\linewidth]{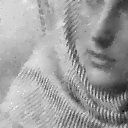}
\includegraphics[width=0.24\linewidth]{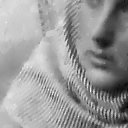}
\includegraphics[width=0.24\linewidth]{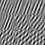}

\caption{\label{fig:denoising_results} Denoising and atom-learning. From left to right: Noisy data, TGV-based reconstruction, proposed method (only \txt{} for the fist image),  nine most important learned atoms.}

\end{figure}

\subsection{Deconvolution} \label{sec:deconvolution}
This section deals with the learning of image features and image reconstruction an an inverse problem setting, where the forward operator is given as a convolution with a Gaussian kernel (standard deviation 0.2, kernel size $9\times 9$ pixels) and the data is degraded by Gaussian noise with standard deviation 0.05 times the image range.  The data fidelity in this case is
\[ D_1(u,f_0)=  \frac{1}{2}\|Au-f_0\|_2^2, \qquad D_2\equiv 0,\]
with $A$ being the convolution operator.

We show results for the \emph{Mix} and the \emph{Barbara} image and compare to \tgv. It can be seen that the improvement is comparable to the denoising case. In particular, the method is still able to learn reasonable atoms from the given, blurry data and in particular for the texture parts the improvement is quite significant. Regarding the optimal PSNR values, the improvement is roughly around 1 to 1.5 decibel. 

\begin{figure}[t]
\center

\includegraphics[width=0.24\linewidth]{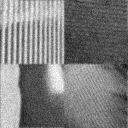}
\includegraphics[width=0.24\linewidth]{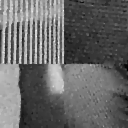}
\includegraphics[width=0.24\linewidth]{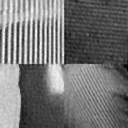}
\includegraphics[width=0.24\linewidth]{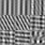}

\includegraphics[width=0.24\linewidth]{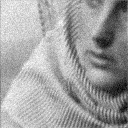}
\includegraphics[width=0.24\linewidth]{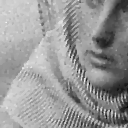}
\includegraphics[width=0.24\linewidth]{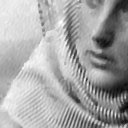}
\includegraphics[width=0.24\linewidth]{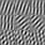}

\caption{\label{fig:recon} Reconstruction from blurry and noisy data. From left to right: Data, TGV, proposed, learned atoms.}

\end{figure}

\subsection{Comparison} \label{sec:comparison}

\begin{figure}[t]
\center 

\begin{tabularx}{\textwidth}{X@{\hspace*{0.1cm}}X@{\hspace*{0.1cm}}X@{\hspace*{0.1cm}}X}

\includegraphics[width=\linewidth]{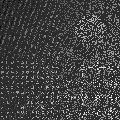} &
\includegraphics[width=\linewidth]{images_l2/img_inpaint/patchtest_recon.png} &
\includegraphics[width=\linewidth]{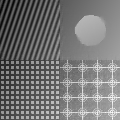} &
\vspace*{-\linewidth}\includegraphics[width=0.99\linewidth]{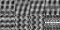}
\includegraphics[width=0.99\linewidth]{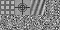} \\

\includegraphics[width=\linewidth]{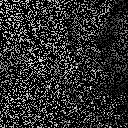} &
\includegraphics[width=\linewidth]{images_l2/img_inpaint/barbara_crop_recon.png} &
\includegraphics[width=\linewidth]{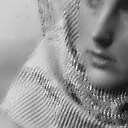} &
\vspace*{-\linewidth}\includegraphics[width=0.99\linewidth]{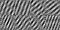}
\includegraphics[width=0.99\linewidth]{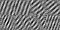}

\end{tabularx}

\caption{\label{fig:semiconv_inp} Comparison of \ctc{} and \ctsc{} for inpainting with 30\% of the pixels given. From left to right: Data, convex, semi-convex, convex atoms (top), semi-convex atoms (bottom).}

\end{figure}

This section compares the method \ctc{} to its semi-convex variant \ctsc{} and to other methods. At first, we consider the learning of atoms from incomplete data and image inpainting in Figure \ref{fig:semiconv_inp}. It can be seen there that for the \emph{Patches} image, the semi-convex variant achieves an almost perfect results: It is able to learn exactly the three atoms that compose the texture part of the image and to inpaint the image very well. For the \emph{Barabara} image, where more atoms are necessary to synthesize the texture part, the two methods yield similar results and also the atoms are similar. These results are also reflected in the PSNR values of Table \ref{tbl:psnr_results}, where \ctsc{} is more that 7 decibel better for the \emph{Patches} image and achieve only a slight improvement for \emph{Barbara}.

Next we consider the semi-convex variant \ctsc{} for denoising the \emph{Patches} and \emph{Barbara} images of Figure \ref{fig:denoising_results}. In this setting, also other methods are applicable and we compare to a costume implementation of a variant of the convolutional lasso algorithm (called \cl) and to BM3D denoising \cite{dabov2007bm3d} (called \bmtd). For the former, we strive to solve the non-convex optimization problem
\[ \min_{ u, (c_i)_i, (p_i)_i } \TV_\rho(u- \sum_{i=1}^k c_i * p_i ) +  \sum_{i=1}^k \|c_i\|_1 + \left\| u - f_0 \right \|_2 ^2  \qquad \text{s.t. } \|p_i\| _2 \leq 1, \int p_i = 0 \text{ for } i=1,\ldots,k\]
where $(c_i)_i$ are coefficient images, $p_i$ are atoms and $k$ is the number of used atoms. Note that we use the same boundary extension, atom-domain-size and stride variable than in the methods \ctc{}, \ctsc{}, and that $\TV_\rho$ denotes a discrete TV functional with a slight smoothing of the $L^1$ norm to make it differentiable (see the source code \cite{convex_lifting_code} for details). For the solution we use an adaption of the algorithm of \cite{pock16ipalm_mh}. For BM3D we use the implementation obtained from \cite{bm3dcode}.

\begin{rem}
We note that, while we provide the comparison to BM3D in order to have a reference on achievable denoising quality, we do not aim to propose an improved denoising method that is comparable to BM3D. In contrast to BM3D, our method constitutes a variational (convex) approach, that is generally applicable for inverse problems and for which we were able to provide a detailed analysis in function space such that in particular stability and convergence results for vanishing noise can be proven. Furthermore, beyond mere image reconstruction, we regard the ability of simultaneous image-atom-learning and cartoon-texture decomposition as an important feature of our approach.
\end{rem}

\begin{figure}[t]
\center 

\begin{tabularx}{\textwidth}{X@{\hspace*{0.1cm}}X@{\hspace*{0.1cm}}X@{\hspace*{0.1cm}}X}

\includegraphics[width=\linewidth]{images_l2/img_denoise/patchtest_recon.png} &
\includegraphics[width=\linewidth]{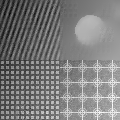} &
\includegraphics[width=\linewidth]{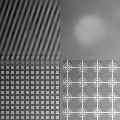} &
\includegraphics[width=\linewidth]{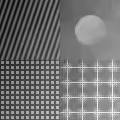}  \\

\includegraphics[width=\linewidth]{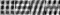} &
\includegraphics[width=\linewidth]{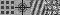} &
\includegraphics[width=\linewidth]{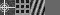} & \\

\includegraphics[width=\linewidth]{images_l2/img_denoise/barbara_crop_recon.png} &
\includegraphics[width=\linewidth]{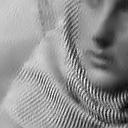} &
\includegraphics[width=\linewidth]{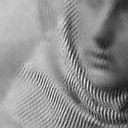} &
\includegraphics[width=\linewidth]{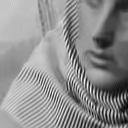} \\

\includegraphics[width=\linewidth]{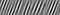} &
\includegraphics[width=\linewidth]{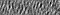} &
\includegraphics[width=\linewidth]{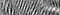} &

\end{tabularx}

\caption{\label{fig:comparison_denoising} Comparison of different methods for denoising the \emph{Patches} and \emph{Barbara} images from Figure \ref{fig:denoising_results}. From left to right: \ctc{}, \ctsc{}, \cl{} and \bmtd{}. The four most important learned atoms are shown below the images.}

\end{figure}

Results for the \emph{Patches} and \emph{Barbara} image can be found in Figure \ref{fig:comparison_denoising}, where for \cl{} we allowed for three atoms for the \emph{Patches} images an tested 3,5 and 7 atoms for the \emph{Barbara} image, showing the best result that was obtained with 7 atoms. It can be seen that, as with the inpainting results, \ctsc{} achieves a very strong improvement compared to \ctc{} for the \emph{Patches} image (obtaining the atoms almost perfectly) and only a slight improvement for the \emph{Barbara} image. The \cl{} method performs similar but slightly worse than \ctsc{}. While for the \emph{Patches} image also the three main features are identified correctly, they are not centered which leads to artifacts in the reconstruction and might be explained by the method being stuck in a local minimum. For the \emph{Patches} image, the result of \bmtd{} are comparable but slightly smoother than the ones of \ctsc{}. In particular, the target-like structure in the bottom left is not very well reconstructed with \bmtd{} but suffers from less remaining noise. For the \emph{Barbara} image, \bmtd{} delivers the best and result, but a slight overs-smoothing is visible. Regarding the PSNR values of Table \ref{tbl:psnr_results}, \bmtd{} performs best and \ctsc{} second best (better that \cl{}), where in accordance with the visual results the difference of \bmtd{} and \ctsc{} is not as high as with \emph{Barbara}.

\begin{figure} 
\begin{tabularx}{\textwidth}{X>{\centering\arraybackslash}p{2cm}>{\centering\arraybackslash}p{2cm}>{\centering\arraybackslash}p{2cm}>{\centering\arraybackslash}p{2cm}}

 & \textbf{Texture} & \textbf{Patches} & \textbf{Mix} & \textbf{Barbara} \\ \toprule 
 \textbf{Inpainting} 
 \\ \midrule 
\tgv   & 10.32 & 19.28 & 20.19 &20.58 \\
\txt / \ctc & \textbf{17.59} & 25.55 & \textbf{23.38} & 23.48 \\
\ctsc   & & \textbf{32.74} & & \textbf{23.6}  \\ \midrule
\textbf{Denoising} \\ \midrule
\tgv  & 11.83 & 23.96 & 23.74 &23.99 \\
\txt / \ctc &  \textbf{16.06} & 25.91 & \textbf{26.07} & 25.0 \\
\ctsc  & & 29.4 & & 25.56 \\
\cl  & & 29.09 & & 25.14 \\
\bmtd  & & \textbf{30.82} & & \textbf{28.15}\\  \midrule
\textbf{Deconvolution} \\ \midrule
\tgv  & & & 23.58 & 22.86\\
\ctc  & & & \textbf{25.19} & \textbf{23.93} \\ \bottomrule

\end{tabularx}
\caption{\label{tbl:psnr_results} Best PSNR result achieved with each method for the parameter test range as specified in Table \ref{tbl:parameter_choice}. The best achieve result for each experiment is written in bold.}
\end{figure}
{
\setlength{\tabcolsep}{3.65pt}
\begin{figure}
\begin{tabular}{@{}rcccccccccccccccccc@{}}
\toprule
& \multicolumn{3}{c}{\ctc} & \phantom{a}
& \multicolumn{4}{c}{\ctsc} & \phantom{a} 
& \multicolumn{1}{c}{\tgv} & \phantom{a} 
& \multicolumn{2}{c}{\txt} & \phantom{a} 
& \multicolumn{1}{c}{\bmtd} & \phantom{a} 
& \multicolumn{2}{c}{\cl}   \\

\cmidrule{2-4} \cmidrule{6-9} \cmidrule{11-11} \cmidrule{13-14}  \cmidrule{16-16} \cmidrule{18-19} \\

 &  $\lambda$ & $\mu $& $\nu $ & & $\lambda$ & $\mu$ & $\nu$ &  $\epsilon$ &  & $\lambda$ &  & $\lambda$ & $\nu$  & & $\lambda $ & & $\lambda$ & $\mu$
\\ \midrule

Decomp. & -  & opt & 0.95 & & & & & & & & & - & 0.75  \\
Inp. & - & opt & 0.975 & &
- & opt & 0.975 & 0.1 & &
-  & &
- & 0.975  \\ 
Den. & opt & opt & 0.975 & &
opt &  opt & 0.975 &  2.0 & & 
 opt & &
opt & 0.975 &  &
opt &&
opt & opt  \\
Deconv. & opt & opt & 0.975 & & & & & & & opt & & & 
\\ \bottomrule
\end{tabular}
\caption{\label{tbl:parameter_choice}Parameter choice for all methods and experiments used in the paper. Here, $\lambda$ always defines the trade-off between data fidelity and regularization, $\mu$ defined the trade-off between cartoon and texture, $\nu$ defined the trade-off between the 1/2 norm and the penalization of singular values and $\epsilon$ defines the degree of non-convexity for the semi-convex potential. Whenever a parameter was optimized over a certain range for each experiment, we write \emph{opt}.}
\end{figure}
}

\section{Discussion}
Using lifting techniques, we have introduced a (potentially convex) variational approach for learning image atoms from corrupted and/or incomplete data. An important part of our work is the analysis of the proposed model, which shows well-posedness results for the proposed model in function space for a general inverse problem setting. The numerical part shows that indeed our model can effectively learn image atoms from different types of data. While this works well also in a convex setting, moving to a semi-convex setting (which is also captured by our theory) yields a further, significant improvement. While the proposed method can also be regarded solely as image reconstruction method, we believe its main feature is in fact the ability to learn image atoms from incomplete data in a mathematically well understood framework.

Interesting future research questions are for instance the exploration of our method for classification problems or the exploration of similar lifting techniques for a mathematical understanding of deep neural networks.

\section{Acknowledgements} MH acknowledges support by the Austrian Science Fund (FWF) (Grant J 4112). TP is supported by the European Research Council under the Horizon 2020 program, ERC starting grant agreement 640156.
\appendix
\section{Appendix: Tensor spaces}

We recall here some basic results on tensor products of Banach spaces that will be relevant for our work. Most of these results are obtained from \cite{Ryan,DiestelUhl}, to which we refer to for further information and a more complete introduction to the topic. 

Throughout this section, let always $(X,\|\cdot \|_X),(Y,\|\cdot \|_Y), (Z,\|\cdot \|_Z)$ be Banach spaces. By $X^*$ we denote the analytic dual of $X$, i.e., the space of bounded linear functionals from $X$ to $\R$. By $\Lin(X,Y)$ and $\B(X\times Y,Z)$ we denote the spaces of bounded linear and bilinear mappings, respectively, where the norm for the latter is given by $\|B\|_\B = \sup \{ \|B(x,y)\|_Z \st \|x\|_X \leq 1,\, \|y\|_Y \leq 1 \}$. In case the image space is the reals, we write $\Lin(X)$ and $\B(X\times Y)$.

\noindent \textbf{Algebraic tensor product.} The tensor product $x \otimes y$ of two elements $x\in X$, $y \in Y$ can be defined as a linear mapping on the space of bilinear forms on $X \times Y$ via
\[ x \otimes y (A) = A(x,y) \]
The algebraic tensor product $X \otimes Y$ is then defined as the subspace of the space of linear functionals on $B(X,Y)$ spanned by elements $x \otimes y$ with $x\in X$, $y \in Y$. 

\noindent \textbf{Tensor norms.} We will use two different tensor norms, the \emph{projective} and the \emph{injective} tensor norm (also known as the largest and smallest reasonable cross norm, respectively). The projective tensor norm on $X \otimes Y$ is defined for $u \in X \otimes Y$ as 
\[ \|u\|_\pi := \inf \left\{ \sum_{i=1}^n  \|x_i\|_X \|y_i\|_Y \st u = \sum_{i=1}^n x_i \otimes y_i , \, n \in \N\right\} .\]
Note that indeed $\|\cdot \|_\pi$ is a norm and $\|x \otimes y\|_\pi = \|x\|_X \|y\|_Y$ (see \cite[Proposition 2.1]{Ryan}). We denote by $X \ctp Y$ the completion of the space $X \otimes Y$ equipped with this norm.
The following result gives a useful representation of elements in $X \ctp Y$ and their projective norm.
\begin{prop} For $u  \in X \ctp Y$ and $\epsilon > 0$ there exist bounded sequences $(x_n)_n \subset X$, $(y_n)_n \subset Y$ such that
\[ u = \sum_{i=1}^\infty x_n \otimes y_n \quad \text{and} \quad \sum_{i=1}^\infty \|x_n\|_X\|y_n\|_Y < \|u\|_\pi + \epsilon. \]
In particular,
\[ \|u\|_\pi = \inf \left \{ \sum_{i=1}^\infty \|x_i\|_X\|y_i\|_Y \st u = \sum_{i=1}^\infty x_i \otimes y_i \right \} .\]
\end{prop}

Now for the injective tensor norm, we note that elements of the tensor product $X \otimes Y$ can be viewed as bounded bilinear forms on $X^* \times Y^*$ by associating to a tensor $u= \sum_{i=1}^n x_i \otimes y_i$ the bilinear form $B_u(\phi,\psi)= \sum_{i=1}^n \phi(x_i)\psi(y_i)$, where this association is unique (see \cite[Section 1.3]{Ryan}). 
Hence $X \otimes Y$ can be regarded as a subspace of $\B(X^* \times Y^*)$ and the injective tensor norm is the norm induced by this space. Thus for $u= \sum_{i=1}^n x_i \otimes y_i$ the injective tensor norm $\|\cdot \|_\itn$ is given as
\[ \|u\|_\itn=  \sup \left\{ \left| \sum_{i=1}^n \phi(x_i)\psi(y_i) \right| \st \|\phi\|_{X^*} \leq 1, \, \|\psi \|_{Y^*}  \leq 1\right\}  \]
and the injective tensor product $X \itp Y$ is defined as the completion of $X \otimes Y$ with respect to this norm.

\noindent \textbf{Tensor lifting.} The next result (see \cite[Theorem 2.9]{Ryan}) shows that there is a one-to-one correspondence between bounded bilinear mappings from $X \times Y$ to $Z$ and bounded linear mappings from $X \ctp Y$ to $Z$.
\begin{prop} For  $B \in \B(X \times Y ,Z)$ there exists a unique linear mapping $\hat{B}:X \ctp Y \rightarrow Z$ such that $\hat{B}(x \otimes y) = B(x,y)$. Further $\hat{B}$ is bounded and the mapping $B \mapsto \hat{B} $ is an isometric isomorphism between the Banach spaces $\B(X\times Y,Z)$ and $\Lin(X \ctp Y,Z)$.
\end{prop}
Using this isometry, for $B \in \B(X\times Y,Z)$ we will always denote by $\hat{B}$ the corresponding linear mapping on the tensor product.

The following result is provided in \cite[Proposition 2.3]{Ryan} and deals with the extension of linear operators to the tensor product.
\begin{prop} Let $S \in \Lin(X,W)$, $T \in \Lin(Y,Z)$. Then there exists a unique operator $S \tp T: X \ctp Y \rightarrow W \ctp Z$ such that $S \tp T(x \otimes y) = (Sx) \otimes (Ty)$. Furthermore, $\|S \tp T \| = \|S\| \|T\|$. 
\end{prop}

\noindent \textbf{Tensor space isometries.}
The following proposition deals with duality of the injective and the projective tensor products. To this aim, we need the notion of \emph{Radon Nikod\'ym property} and \emph{approximation property}, which we will not define here but rather refer to \cite[Sections 4 and 5]{Ryan} and \cite{DiestelUhl}. For our purposes, it is important to note that both properties hold for $L^r$-spaces with $r \in (1,\infty)$, the \emph{Radon Nikod\'ym property} holds for reflexive spaces, but while we cannot expect the \emph{Radon Nikod\'ym property} to hold for $L^\infty$ and $ \M$, the approximation property does.
\begin{lem}\label{lem:dual_of_itp_is_ptp} Assume that either $X^*$ or $Y^*$ has the Radon Nikod\'ym property and that either $X^*$ or $Y^*$ has the approximation property. Then
\[  (X \itp Y)^* \hat{=} X^* \ctp Y^*\]
and for simple tensors $u = \sum_{i=1}^n x_i \otimes y_i \in X \itp Y$ and $u^* = \sum_{i=1}^m x^*_i \otimes y^*_i \in X^* \ctp Y^*$ the duality pairing is given as
\[ \langle u^*,u\rangle  = \sum_{i=1}^n \sum_{j=1}^m \langle x_j^*,x_i \rangle \langle y_j^*,y_i \rangle \]
\proof
The identification of the duals is shown in \cite[Theorem 5.33]{Ryan}. For the duality paring, we first note that the action of an element $u^* \in X^* \ctp Y^*$ on $X \itp Y$ is given as the action of the associated bilinear form $B_{u^*}$ \cite[Section 3.4]{Ryan}, which for simple tensors $u= \sum_{i=1}^n x_i \otimes y_i$ can be given as
\[ \langle B_{u^*},u\rangle = \sum_{i=1}^n B_{u^*} (x_i,y_i) .\]
Now in case also $u^*$ is a simple tensor, i.e., $u^* = \sum_{i=1}^m x_i^* \otimes y_i^*$, the action of this bilinear form can be given more explicitly \cite[Section 1.3]{Ryan}, which yields
\[ \langle B_{u^*},u\rangle = \sum_{i=1}^n  \sum_{i=1}^m \langle x^*_i,x_i\rangle \langle y^*_i,y_i\rangle .\] \qedhere
\end{lem}
The duality between the injective and projective tensor product will be used for compactness assertions on subsets of the latter. To this aim, we note in the following lemma that separability of the individual space transfers to the tensor product. As consequence, in case $X^*$ and $Y^ *$ satisfy the assumption of Lemma \ref{lem:dual_of_itp_is_ptp} above and both admit a separable predual, also $X^ * \ctp Y^ *$ admits a separable predual and hence bounded sets are weakly* compact.
\begin{lem} \label{lem:tensor_product_separable} Let $X,Y$ be separable spaces. Then both $X \itp Y$ and $X \ctp Y$ are separable.
\proof 
Take $X'$ and $Y'$ to be dense countable subsets of $X$ and $Y$, respectively. First note that it suffices to show that any simple tensor $x \otimes y$ can be approximated arbitrarily close by $x'\otimes y'$ with $x' \in X' $, $y' \in Y'$. But this is true since (using \cite[Propositions 2.1 and 3.1]{Ryan})
\[ \| x \otimes y - x'\otimes y'\| \leq \| x \otimes y - x\otimes y'\|  + \| x \otimes y' - x'\otimes y'\| = \|x\|\|y-y'\| + \|y'\|\|x-x'\| ,\]
where $\|\cdot \|$ denotes either the projective or the injective norm.

\end{lem}
The following result, which can be obtained by direct modification of the result shown at the beginning of \cite[Section 3.2]{Ryan}, provides an equivalent representation of the injective tensor product in a particular case.
\begin{lem} \label{lem:isometriy_itp_cinfty} Denote by $C_c(\Omega_\Sigma,X)$ the space of compactly supported continuous functions mapping from $\Omega_\Sigma$ to $X$ and denote by $C_0(\Omega_\Sigma,X)$ its completion with respect to the norm $\|\phi\|_\infty := \sup_{t \in \Omega_\Sigma} \|X\|_X$. Then we have that
\[  C_0(\Omega_\Sigma) \itp X  \hat{=} C_0(\Omega_\Sigma,X) \]
where the isometry is given as the completion of the isometric mapping $J: C_0(\Omega_\Sigma) \otimes X \rightarrow C_0(\Omega_\Sigma,X)$ defined for $u = \sum_{i=1}^n f_i \otimes x_i$ as
\[ Ju(t):= \sum _{i=1}^n f_i(t)x_i .\]
\end{lem}

Next we consider the identification of tensor products with linear operators which is provided in the following proposition \cite[Corollary 4.8]{Ryan}.
\begin{prop} \label{prop:nuclear_norm_operator} Define the mapping $J:X^* \ctp Y \rightarrow \Lin (X,Y)$ as
\[ u= \sum_{i=1}^\infty \phi_n \otimes y_n \mapsto L_u:X \rightarrow Y \text{ where } L_u(x)= \sum_{i=1}^\infty \phi_n(x) y_n .\]
Then, $J$ is well-defined and has unit norm. Defining $\Nuc(X,Y)\subset \Lin(X,Y)$ as the range of $J$, equipped with the norm
\[ \|T\|_\nuc = \inf \bigg\{ \sum_{i=1}^\infty \|\phi_n\|_{X^*}\|y_n\|_Y \st T(x) = \sum_{i=1}^\infty \phi_n(x)y_n \bigg \}, \]
we get that $\Nuc(X,Y)$ is a Banach space, called the space of nuclear operators. If further either $X^*$ or $Y$ has the approximation property, then $J$ is an isometric isomorphism, that is, we can identify
\[ X^* \ctp Y = \Nuc(X,Y) \]
\end{prop}
It is easy to see that nuclear operators are compact that that we can equivalently write
\[  \|T\|_\nuc = \inf \bigg\{ \sum_{i=1}^\infty \sigma_i \st T(x) = \sum_{i=1}^\infty \sigma_i \phi_i(x)y_i, \, \|\phi_i\|_{X^*}\leq 1,\, \|y_i\|_Y\leq 1 \bigg \}. \]
Also, in a Hilbert space setting (see \cite{Weidmann2012linear} for details), it
is a classical result that for any compact $T \in \Lin (H_1,H_2)$ with $(H_1,(\cdot,\cdot)),(H_2(\cdot,\cdot))$ Hilbert spaces there exist orthonormal systems $(x_i)_i$, $(y_i)_i$ and uniquely defined singular values $(\sigma_i)_i := (\sigma_i(T))_i$ such that \[ Tx = \sum_{i=1}^\infty \sigma_i (x,x_i)y_i .\]
In addition, in case $T$ has finite nuclear norm, it follows that $ \|T\|_\nuc = \sum_{i=1}^\infty \sigma _i .$

\bibliography{lit_dat}
\bibliographystyle{abbrv}

\end{document}